\documentclass[a4paper,11pt]{article}
\usepackage[noabacus]{mfpaperstuff2}
\begin{document}
\renewcommand\labelenumi{(\theenumi)}
\newcommand\wst{\weyl s\times\weyl t}
\newcommand\sset[1]{\calq(#1)}
\newcommand\stab{\operatorname{Stab}}
\newcommand\jn[2]{#1\star#2}
\newcommand\xx[1]{\card{#1}_{s,t}}
\newcommand\ci[2]{#1{\circ}#2}
\newcommand\dst[4]{#1\stackrel{#3,#4}\leftrightharpoons#2}
\newcommand\cst[4]{#1\stackrel{#3,#4}\rightharpoondown#2}
\newcommand\imps[2]{\item[{\rm(\ref{#1}$\Rightarrow$\ref{#2})}]}
\newcommand\impss[2]{\item[{\rm(\ref{#1}$\Leftrightarrow$\ref{#2})}]}
\newcommand\ol[1]{\widebar{#1}}
\newcommand\pvt{peak\xspace}
\newcommand\pvts{peaks\xspace}
\newcommand\bkn[3]{\shade[shading=ball,ball color=white](#1,#2)circle(.4);\draw(#1,#2);\draw(#1,#2)node{$\vphantom{1}\smash{#3}$};}
\newcommand\whn[3]{\draw(#1,#2)node{$\vphantom{1}\smash{#3}$};}
\newcommand\bn[2]{\shade[shading=ball,ball color=white](#1,#2)circle(.4);}
\newcommand\wn[2]{\draw(#1,#2)node{.};}
\newcommand\rsx[1]{\hat{#1}}
\newcommand\zom{$(0,1)$-matrix\xspace}
\newcommand\zoms{$(0,1)$-matrices\xspace}
\newcommand\dg[2]{g_{#1}(#2)}
\newcommand\zsz{\bbz/s\bbz}
\newcommand\ztz{\bbz/t\bbz}
\newcommand\weyl[1]{\calw_{#1}}
\newcommand\hyp[1]{\calh_{#1}}
\newcommand\cor[2]{\operatorname{c}_{#1}(#2)}
\newcommand\wei[2]{\operatorname{w}_{#1}(#2)}
\renewcommand\be[1]{\calb(#1)}
\renewcommand\sc{self-conjugate\xspace}
\newcommand\scp{\sc partition\xspace}
\newcommand\scps{\scp{}s\xspace}
\newcommand\pst[3]{\calp_{#2,#3}(#1)}
\newcommand\pstn{\pst n\si \tau}
\newcommand\stm{$(s,t)$-minimal\xspace}
\newcommand\stmp{\stm partition\xspace}
\newcommand\stmps{\stmp{}s\xspace}
\newcommand\sitam{$(\si,\tau)$-minimal\xspace}
\newcommand\sitamp{\sitam partition\xspace}
\newcommand\sitamps{\sitamp{}s\xspace}
\newcommand\pp{pinch-point\xspace}
\newcommand\pps{\pp{}s\xspace}
\newcommand\mn[2]{\mathcal{M}_{#1,#2}}
\newcommand\mns{\mn\si\tau}
\newcommand\nn[2]{\mathtt{m}_{#1,#2}}
\newcommand\nns{\nn\si\tau}
\newcommand\cores[1]{\calc_{#1}}

\title{Minimal partitions with a given $s$-core and $t$-core}
\msc{05A17, 05E10, 05E18}

\toptitle

\begin{abstract}
Suppose $s$ and $t$ are coprime positive integers, and let $\si$ be an $s$-core partition and $\tau$ a $t$-core partition. In this paper we consider the set $\pstn$ of partitions of $n$ with $s$-core~$\si$ and $t$-core~$\tau$. We find the smallest $n$ for which this set is non-empty, and show that for this value of $n$ the partitions in $\pstn$ (which we call \emph{\sitam} partitions) are in bijection with a certain class of \zoms with $s$ rows and $t$ columns.

We then use these results in considering conjugate partitions: we determine exactly when the set $\pstn$ consists of a conjugate pair of partitions, and when $\pstn$ contains a unique \scp.
\end{abstract}

\section{Introduction}

In recent years there has been considerable interest in the study of \emph{core partitions}. For any positive integer $s$, an integer partition is an \emph{$s$-core partition} (or simply an \emph{$s$-core} in this paper) if it has no hooks of length $s$. In the case where $s$ is a prime, $s$-cores play an important role in the representation theory of the symmetric groups in characteristic $s$. These partitions are also important in number theory, and more recently have been studied from the point of view of enumerative combinatorics. A popular theme is to take a second positive integer $t$, and study the set of \emph{$(s,t)$-cores}, i.e.\ partitions which are both $s$- and $t$-cores. More generally, one can consider the set of partitions which are $s$-cores for all $s$ in some given set; if the elements of this set are coprime, there are only finitely many such partitions, leading to various enumerative results. See \cite{and,trip,ahj,johns} for some of these results; there are many more results in which further restrictions are placed on the partitions.

Given a partition $\la$ which is not an $s$-core, one can define the $s$-core of $\la$ by repeatedly removing rim hooks of length $s$. This is also significant in the representation theory of the symmetric group in characteristic $s$: by the Brauer--Robinson theorem \cite{brauer,rob}, the $s$-core of a partition determines the $s$-block in which the corresponding irreducible representation lies. The focus in the present paper is to study the $s$-core and $t$-core operations on partitions simultaneously. This approach has been taken before, in \cite{olss,nath,mfcores,mfgencores}. In particular, in \cite{mfgencores} a family of partitions called \emph{$[s{:}t]$-cores} was introduced. One of the equivalent characterisations of these partitions is that a partition $\la$ is an $[s{:}t]$-core \iff there is no other partition with the same size, $s$-core and $t$-core as $\la$; this implies also that there is no smaller partition with the same $s$-core and $t$-core as $\la$. In the present paper we pursue this idea further: we take an $s$-core $\si$ and a $t$-core $\tau$, and study the set $\pstn$ of partitions of $n$ with $s$-core $\si$ and $t$-core $\tau$. We determine the smallest $n$ for which $\pstn$ is non-empty, and show that for this value of $n$ the set $\pstn$ (whose elements we call \emph{\sitam} partitions) is in bijection with a class of \zoms. We use this to recover the result from \cite{mfgencores} determining exactly when $\card{\pstn}=1$, and in addition we determine when $\card{\pstn}=2$.

We then use these results in relation to conjugation of partitions. The \emph{conjugate} of a partition $\la$ is the partition obtained by reflecting the Young diagram on the diagonal, and conjugation is significant in the representation theory of the symmetric and alternating groups. We prove two results analogous to our result from \cite{mfgencores}: we determine for which $n$ there is a unique partition of $n$ up to conjugation with $s$-core $\si$ and $t$-core $\tau$, and we determine for which $n$ there is a unique \scp of $n$ with $s$-core $\si$ and $t$-core $\tau$.

In proving these results, we concentrate initially on the case where $\si$ and $\tau$ are both $(s,t)$-cores. We then consider actions of affine symmetric and hyperoctahedral groups on the set of partitions. If we define a partition to be \emph{\stm} if it is \sitam for some $\si,\tau$, then these actions preserve the set of \stmps, which allows us to extend our results to the case where $\si$ and $\tau$ are not necessarily $(s,t)$-cores.

\begin{ack}
The results presented in this paper were obtained with the help of extensive calculations performed in GAP \cite{gap}.
\end{ack}

\section{Partitions, cores and beta-sets}

We begin by introducing various concepts relating to partitions.

\subsection*{Elementary notation}

In this paper, $\bbn$ denotes the set of positive integers, $-\bbn$ the set of negative integers, and $\bbn_0$ the set of non-negative integers.

If $s\in\bbn$, then $\bbz/s\bbz$ is the set of cosets $a+s\bbz$ for $a\in\bbz$. (We do not employ the standard abuse of notation in which $\zsz$ is identified with $\{0,\dots,s-1\}$.) Given any tuple of objects $\ltup{x_i}{i\in\zsz}$ indexed by $\zsz$, we may write $x_a$ to mean $x_{a+s\bbz}$ for any integer $a$.

If $X$ is any finite set of integers, we write $\sum X$ for the sum of the elements of $X$.

\subsection*{Partitions}

A \emph{composition} is an infinite sequence $\la=(\la_1,\la_2,\dots)$ of non-negative integers with finite sum. We write $\card\la=\la_1+\la_2+\dots$, and we say that $\la$ is a composition of $\card\la$. We say that $\la$ is a \emph{partition} if $\la_1\gs\la_2\gs\dots$. We write $\calp$ for the set of all partitions, and $\calp(n)$ for the set of all partitions of $n$.

When writing a partition, we usually group together equal parts and omit trailing zeroes, so that the partition $(7,7,7,6,6,6,3,0,0\dots)$ is written as $(7^3,6^3,3)$. We write $\varnothing$ for the unique partition of $0$. The \emph{Young diagram} of a partition $\la$ is the set
\[
\lset{(r,c)\in\bbn^2}{c\ls\la_r},
\]
whose elements we call the \emph{nodes} of $\la$. We draw Young diagrams using the English convention, in which the partition $(4^2,1)$ is depicted as follows.
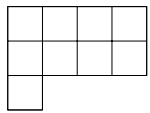
\[
\yng(4^2,1)
\]

If $\la$ is a partition, the \emph{conjugate} partition is the partition $\la'$ obtained by reflecting the Young diagram of $\la$ on the main diagonal. In other words, the partition given by
\[
\la'_r=\cardx{\lset{c\in\bbn}{\la_c\gs r}}.
\]
We say that $\la$ is \emph{\sc} if $\la'=\la$.

The \emph{dominance order} is a partial order defined on the set of all partitions by writing $\la\dom\mu$ (and saying that $\la$ \emph{dominates} $\mu$) if $\card\la=\card\mu$ and $\la_1+\dots+\la_r\gs\mu_1+\dots+\mu_r$ for all $r$.

Now we define the \emph{beta-set} of a partition $\la$. This is the set
\[
\be\la=\lset{\la_r-r}{r\in\bbn}.
\]
The next \lcnamecref{basicbeta} (which is easy to prove by induction) gives some basic information about beta-sets.

\begin{lemma}\label{basicbeta}
Suppose $\la\in\calp$. Then $\bbn_0\cap\be\la$ and $(-\bbn)\setminus\be\la$ are finite sets of the same size, and
\[
\card\la=\sum(\bbn_0\cap\be\la)-\sum((-\bbn)\setminus\be\la).
\]
Conversely, any set $B\subset\bbz$ for which $\bbn_0\cap B$ and $(-\bbn)\setminus B$ are finite sets of the same size is the beta-set of a unique partition.
\end{lemma}

Conjugation of partitions is also encoded in beta-sets: it is easy to prove that
\[
\be{\la'}=\bbz\setminus\lset{-1-b}{b\in\be\la}
\]
for any partition $\la$.

\subsection*{Rim hooks and cores}

The \emph{rim} of a partition $\la$ is the set of nodes $(r,c)$ of $\la$ such that $(r+1,c+1)$ is not a node of $\la$. If $s\in\bbn$, a \emph{rim $s$-hook} of $\la$ is a connected set of $s$ nodes in the rim which can be removed to leave the Young diagram of a smaller partition. A partition is an \emph{$s$-core} if it has no rim $s$-hooks. We write $\cores s$ for the set of all $s$-cores.

If $\la$ is any partition, the $s$-core of $\la$ is the partition obtained by recursively removing rim $s$-hooks until none remain. This is defined independently of the choice of rim $s$-hook removed at each stage. We write $\cor s\la$ for the $s$-core of $\la$. We define the \emph{$s$-weight} $\wei s\la$ of $\la$ to be the number of rim $s$-hooks removed to reach the $s$-core. Then $\card\la=\card{\cor s\la}+s\wei s\la$.

Removal of rim hooks is closely related to the beta-set of a partition. It was first shown by Nakayama \cite{naka1} that removing a rim $s$-hook from a partition corresponds to reducing an element of $\be\la$ by $s$. In particular, $\la$ is an $s$-core \iff $b-s\in\be\la$ for all $b\in\be\la$. 

For example, the $5$-core of $(6^2,5,1^2)$ is $(3,1)$, as we see from the following diagrams, which show the successive removal of rim $5$-hooks. We also show the effect of removing these rim hooks on the beta-sets of these partitions (at each stage, the underlined entry is reduced by $5$).
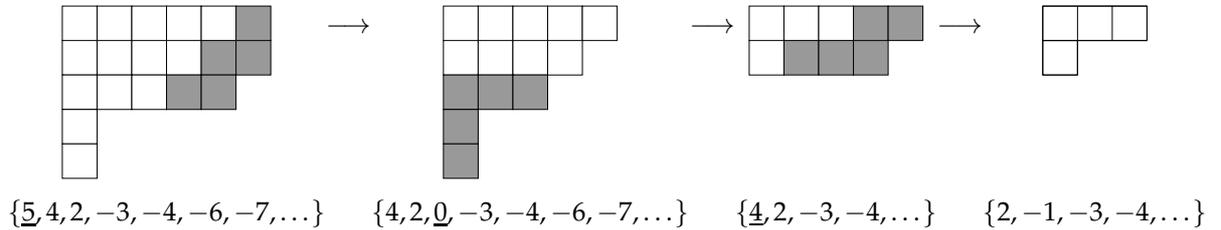
\newcommand\gre{\Yfillcolour{white!60!black}}
\newcommand\whit{\Yfillcolour{white}}
{\small
\[
\begin{array}{c@{}c@{}c@{}c@{}c@{}c@{}c}
\gyoung(;;;;;!\gre;,!\whit;;;;!\gre;;,!\whit;;;!\gre;;,!\whit;,;)&\longrightarrow&\gyoung(;;;;;,;;;;,!\gre;;;,;,;)&\longrightarrow&\gyoung(!\whit;;;!\gre;;,!\whit;!\gre;;;)&\longrightarrow&\yng(3,1)
\\[59pt]
\{\underline5,4,2,-3,-4,-6,-7,\dots\}
&&
\{4,2,\underline0,-3,-4,-6,-7,\dots\}
&&
\{\underline4,2,-3,-4,\dots\}
&&
\{2,-1,-3,-4,\dots\}
\end{array}
\]
}

The connection between rim hooks and beta-sets gives the following statement.

\begin{lemma}\label{corebij}
Suppose $\la,\mu\in\calp$. \Tfae.
\begin{enumerate}
\item\label{smcr}
$\cor s\la=\cor s\mu$.
\item\label{fbij}
There is a bijection $\phi:\be\la\to\be\mu$ such that $\phi(b)\equiv b\ppmod s$ for each $b\in\be\la$, and $\phi(b)=b$ for all but finitely many $b\in\be\la$.
\item\label{sbij}
There is a bijection $\phi:\be\la\setminus\be\mu\to\be\mu\setminus\be\la$ such that $\phi(b)\equiv b\ppmod s$ for each $b\in\be\la\setminus\be\mu$.
\end{enumerate}
\end{lemma}
\begin{pf}\indent
\begin{description}
\vst
\impss{smcr}{fbij}
Let $\si=\cor s\la$ and $\tau=\cor s\mu$. Let's call a bijection $\phi$ as in (\ref{fbij}) a \emph{good} bijection for $\la$ and $\mu$. The fact that $\be\si$ is obtained from $\be\la$ by successively reducing entries by $s$ means that there is a good bijection $\hat\phi$ for $\la$ and $\si$. Similarly there is a good bijection $\check\phi$ for $\mu$ and $\tau$, so if $\si=\tau$ then the composition $\check\phi^{-1}\circ\hat\phi$ is a good bijection for $\la$ and $\mu$.

Conversely, if there is a good bijection $\phi$ for $\la$ and $\mu$, then there is a good bijection $\psi=\check\phi\circ\phi\circ\hat\phi^{-1}$ for $\si$ and $\tau$. For any $i\in\zsz$ the bijection $\psi$ restricts to a bijection $\bar\psi$ from $\be\si\cap i$ to $\be\tau\cap i$ such that $\bar\psi(b)=b$ for all but finitely many elements of $\be\si\cap i$. But the fact that $\si$ and $\tau$ are $s$-cores means that $\be\si\cap i=\{c,c-s,c-2s,\dots\}$ and $\be\tau\cap i=\{d,d-s,d-2s,\dots\}$ for some integers $c,d$. The existence of $\bar\psi$ now implies that $c=d$, so that $\be\si\cap i=\be\tau\cap i$. This applies for all $i$, so that $\be\si=\be\tau$, and hence $\si=\tau$.

\impss{fbij}{sbij}
Given a bijection $\phi:\be\la\setminus\be\mu\to\be\mu\setminus\be\la$ as in (\ref{sbij}), we extend it to a bijection from $\be\la$ to $\be\mu$ by defining $\phi(b)=b$ for all $b\in\be\la\cap\be\mu$. This bijection is good, because $\be\la\setminus\be\mu$ is finite.

Conversely, suppose $\phi$ is a good bijection for $\la$ and $\mu$. Then we claim we can modify $\phi$ so that $\phi(b)=b$ for all $b\in\be\la\cap\be\mu$: if there is some $b\in\be\la\cap\be\mu$ with $\phi(b)\neq b$, then let $a=\phi^{-1}(b)$ and $c=\phi(b)$, and replace $a\mapsto b$ and $b\mapsto c$ in the definition of $\phi$ with $a\mapsto c$ and $b\mapsto b$. This reduces the number of $b\in\be\la\cap\be\mu$ for which $\phi(b)\neq b$, and so in finitely many steps we can reach a bijection $\phi$ in which $\phi(b)=b$ for all $b\in\be\la\cap\be\mu$. Now restricting $\phi$ to $\be\la\setminus\be\mu$ gives a bijection as in (\ref{sbij}).\qedhere
\end{description}
\end{pf}

We also need some basic results on rim hooks and conjugation. These are well known to experts, but it is easier to provide a proof than a reference.

\begin{lemma}\label{basichook}\indent
\begin{enumerate}
\vspace{-\topsep}
\item\label{conjcore}
If $\la\in\calp$, then $\cor s{\la'}=\cor s\la'$. Hence $\la'$ is an $s$-core \iff $\la$ is.
\item\label{scwt1}
If $\si\in\cores s$, then there is at most one \scp with $s$-core $\si$ and $s$-weight $1$.
\end{enumerate}
\end{lemma}

\begin{pfenum}
\item
If $\mu$ is a partition obtained from $\la$ by removing a rim $s$-hook, then by reflecting the diagrams on the diagonal we find that $\mu'$ is obtained from $\la'$ by removing a rim $s$-hook. Now the result follows.
\item
This is most easily seen in terms of beta-sets. If we write $\ol b=-1-b$ for any integer $b$, then (as noted above) a partition $\la$ is \sc \iff $\ol{\be\la}=\bbz\setminus\be\la$.

If $\si\neq\si'$ then by (\ref{conjcore}) there are no \scps with $s$-core $\si$. So assume $\si=\si'$. If $\la$ is a partition with $s$-core $\si$ and $s$-weight $1$, then $\be\la$ is obtained from $\be\si$ by replacing some integer $b$ with $b+s$. Given that $\ol{\be\si}=\bbz\setminus\be\si$, the only way we can also have $\ol{\be\la}=\bbz\setminus\be\la$ is if $b=(-1-s)/2$. So there is only one possibility for~$\la$.
\end{pfenum}

\subsection*{Simultaneous cores}

\begin{framing}{black}
From now on, we fix coprime natural numbers $s$ and $t$ greater than $1$.
\end{framing}

A lot of recent literature has been concerned with studying $s$- and $t$-cores simultaneously. Suppose $\la$ is both an $s$-core and a $t$-core. In this case we say that $\la$ is an \emph{$(s,t)$-core}. We write $\cores{s,t}$ for the set of all $(s,t)$-cores. In fact $\cores{s,t}$ is finite: Anderson \cite[Theorems 1 and 3]{and} showed that $\card{\cores{s,t}}=\mfrac1{s+t}\mbinom{s+t}s$.

$\cores{s,t}$ contains a unique largest partition, denoted $\kappa_{s,t}$. It can be constructed as follows. Let
\[
\calx=\left\{\mfrac{-t-1+s(t-1)}2,\mfrac{-t-1+s(t-3)}2,\dots,\mfrac{-t-1+s(3-t)}2,\mfrac{-t-1+s(1-t)}2\right\}.
\]
Then $\calx$ is a set of $t$ integers which are congruent modulo $s$ but pairwise incongruent modulo $t$, and it is shown in \cite[Section 5]{mfcores} that $\kappa_{s,t}$ is the partition whose beta-set is
\[
\bigcup_{x\in\calx}\{x,x-t,x-2t,\dots\}.
\]
Later we shall need the following result concerning $\kappa_{s,t}$, which appears to be new.

\begin{propn}\label{kappaplusst}
Suppose $\la$ is a partition with $\card\la=\card{\ka_{s,t}}+st$, and that $\cor s\la=\cor t\la=\ka_{s,t}$. Then $\la$ is obtained by adding a rim $st$-hook to $\ka_{s,t}$.
\end{propn}

\begin{pf}
In this proof we write $\ka$ for $\ka_{s,t}$. Let $\calx$ be the set defined above, and for any integer $b$, write $x_b$ for the unique element of $\calx\cap(b+t\bbz)$.

Because $\cor t\la=\ka$, \cref{corebij} says that there is a bijection $\phi:(\be\la\setminus\be\ka)\to(\be\ka\setminus\be\la)$ such that $\phi(b)\equiv b\ppmod t$ for each $b\in\be\la\setminus\be\ka$. Then from \cref{basicbeta}
\[
st=\card\la-\card\ka=\sum_{b\in\be\la\setminus\be\ka}(b-\phi(b))\ =\ \sum_{b\in\be\la\setminus\be\ka}(b-x_b)\ +\ \sum_{b\in\be\ka\setminus\be\la}(x_b-b).\tag*{(\textasteriskcentered)}
\]
The relationship between $\be\ka$ and $\calx$ means that $b>x_b$ for every $b\in\bbz\setminus\be\ka$, while $b\ls x_b$ for every $b\in\be\ka$.

Since $\cor s\la=\ka$, there is another bijection $\psi:(\be\la\setminus\be\ka)\to(\be\ka\setminus\be\la)$ such that $\psi(b)\equiv b\ppmod s$ for all $b$. Because all the elements of $\calx$ are congruent modulo $s$, this means that the integer
\[
y_b=(b-x_b)+(x_{\psi(b)}-\psi(b))
\]
is divisible by $s$ for every $b\in\be\la\setminus\be\ka$. By definition $y_b$ is also divisible by $t$, so it is divisible by $st$. But $b-x_b>0$ and $x_{\psi(b)}-\psi(b)\gs0$ for every $b$, so $y_b$ is a positive multiple of $st$. But (\textasteriskcentered) says that
\[
\sum_{b\in\be\la\setminus\be\ka}y_b=st,
\]
and therefore $\card{\be\la\setminus\be\ka}=1$, which means that $\la$ is obtained from $\ka$ by adding a rim hook.
\end{pf}

When considering $s$- and $t$-cores simultaneously, it is useful to depict beta-sets of partitions using the $(s,t)$-diagram introduced by Anderson \cite{and}. This diagram consists of the integer lattice $\bbz^2$, with the position $(r,c)$ replaced by the integer $rt+cs$. For consistency with Young diagrams, we draw $(s,t)$-diagrams so that the coordinate $r$ increases down the page, and the coordinate $c$ increases from left to right.

The $(s,t)$-diagram of a partition $\la$ is obtained by placing a bead on the diagram at all positions labelled by elements of $\be\la$. For example, a portion of the $(3,4)$-diagram for the partition $(2)$ is as follows.
{\footnotesize
\[
\begin{tikzpicture}[scale=.7,yscale=-1]
\bkn{-1}0{-12}
\bkn00{-9}
\bkn10{-6}
\bkn20{-3}
\whn30{0}
\bkn{-1}1{-8}
\bkn01{-5}
\bkn11{-2}
\bkn211
\whn314
\bkn{-1}2{-4}
\whn02{-1}
\whn122
\whn225
\whn328
\whn{-1}3{0}
\whn03{3}
\whn136
\whn239
\whn33{12}
\end{tikzpicture}
\]
}
Observe that the $(s,t)$-diagram is periodic: it is unchanged by translations through multiples of $(s,-t)$. We can easily tell from the $(s,t)$-diagram whether $\la$ is an $(s,t)$-core: every bead must have a bead immediately above and a bead immediately to the left.

We end this section by recalling Olsson's theorem on cores.

\begin{thm}[\xcite{olss}{Theorem 1}]\label{olssonthm}
Suppose $\si$ is an $s$-core. Then $\cor t\si$ is also an $s$-core.
\end{thm}

\section{Minimal partitions with given cores}\label{stmsec}

\subsection*{Partitions with prescribed cores}

Suppose $\si\in\cores s$ and $\tau\in\cores t$. It is an easy exercise using the Chinese Remainder Theorem to show that there exist partitions with $s$-core $\si$ and $t$-core $\tau$. One of the main aims of this paper is to find the smallest such partition(s). So let's say that a partition with $s$-core $\si$ and $t$-core $\tau$ is \emph{\sitam} if there is no smaller partition with $s$-core $\si$ and $t$-core $\tau$. In general, we say that a partition $\la$ is \emph{\stm} if it is $(\cor s\la,\cor t\la)$-minimal; that is, there is no smaller partition with the same $s$-core and $t$-core as $\la$.

In this section we will attempt to describe the \sitamps for given $\si$ and $\tau$; in particular, to say how big these partitions are, and how many of them there are. Given $\si\in\cores s$ and $\tau\in\cores t$, we write $\mns$ for the set of \sitamps, and $\nns$ for the common size of these partitions. More generally, for any $n$ we write $\pstn$ for the set of partitions of $n$ with $s$-core $\si$ and $t$-core $\tau$.

Determining $\nns$ for given $\si$ and $\tau$ allows us to determine for which $n$ the set $\pstn$ is non-empty. If $\pstn\neq\emptyset$, then by definition $n\gs\nns$. But also $n\equiv\card\si\equiv\nns\ppmod s$ and $n\equiv\card\tau\equiv\nns\ppmod t$, so that $n\equiv\nns\ppmod{st}$. And the conditions $n\gs\nns$ and $n\equiv\nns\ppmod{st}$ are sufficient for the existence of a partition of $n$ with $s$-core $\si$ and $t$-core $\tau$: taking a partition of $\nns$ with $s$-core $\si$ and $t$-core $\tau$, one can add a rim $st$-hook to get a partition $\la^+$ of $\nns+st$. A rim $st$-hook can be decomposed into $s$ rim $t$-hooks (or $t$ rim $s$-hooks) so that $\la^+$ also has $s$-core $\si$ and $t$-core $\tau$. Repeating this operation as many times as needed gives the required partition of $n$. In fact (since there are at least $st$ ways to add a rim $st$-hook to any partition) we deduce the following.

\begin{propn}\label{whenpst}
Suppose $\si\in\cores s$, $\tau\in\cores t$ and $n\in\bbn_0$. Then $\pstn\neq\emptyset$ \iff $n\gs\nns$ and $n\equiv\nns\ppmod{st}$. Furthermore, if $n>\nns$ with $n\equiv\nns\ppmod{st}$, then $\card{\pstn}\gs st$.
\end{propn}

This has the following consequence, which we will need later.

\begin{cory}\label{smallwt}
If $\la\in\calp$ with either $\wei s\la<t$ or $\wei t\la<s$, then $\la$ is \stm.
\end{cory}

\begin{pf}
Let $\si=\cor s\la$ and $\tau=\cor t\la$, and take $\mu\in\mns$. If $\la$ is not \stm, then $\card\la\gs\card\mu+st$ by \cref{whenpst}, so that $\wei s\la\gs\wei s\mu+t\gs t$ and similarly $\wei t\la\gs s$.
\end{pf}

In this section we will mainly concentrate on the case where $\si$ and $\tau$ are both $(s,t)$-cores; in \cref{group1} we will extend some of our results to the general case using group actions.

\subsection*{Rectangles in the $(s,t)$-diagram}

It will help us to introduce certain subsets of $\bbz$. Given $x\in\bbz$, define
\begin{align*}
\calu_x&=\lset{x+as+bt}{a,b>0},\\
\call_x&=\lset{x-as-bt}{a,b\gs0},\\
\calr_x&=\lset{x+as-bt}{1\ls a\ls t,\ 0\ls b\ls s-1}.
\end{align*}
Call the $\calr_x$ the \emph{$x$-rectangle}. Observe that $\bbz$ is the disjoint union of $\calu_x$, $\call_x$ and $\calr_x$, and that $\calr_x$ is a transversal of the congruence classes modulo $st$.

The sets $\calr_x$, $\calu_x$ and $\call_x$ are most easily visualised in the $(s,t)$-diagram: the elements of $\calr_x$ form a rectangle (in fact, a sequence of repetitions of a rectangle, given the periodicity of the diagram). $\calu_x$ is the region below and to the right of these rectangles, and $\call_x$ is the region above and to the left.

For example, in the $(3,4)$-diagram with $x=1$, we have the following picture (with $\calr_1$ shaded).
{\footnotesize
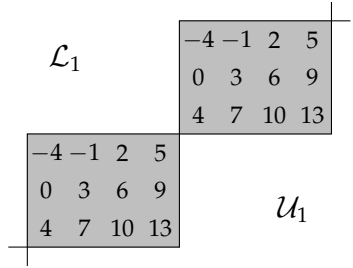
\[
\begin{tikzpicture}[scale=.5]
\filldraw[gray!50!white](-.5,-.5)rectangle++(4,3);
\filldraw[gray!50!white](3.5,2.5)rectangle++(4,3);
\draw(-1,-.5)--++(4.5,0)--++(0,6)--++(4.5,0);
\draw(-.5,-1)--++(0,3.5)--++(8,0)--++(0,3.5);
\draw(0,0)node{$4$}++(1,0)node{$7$}++(1,0)node{$10$}++(1,0)node{$13$}++(0,1)node{$9$}++(-1,0)node{$6$}++(-1,0)node{$3$}++(-1,0)node{$0$}++(0,1)node{$-4$}++(1,0)node{$-1$}++(1,0)node{$2$}++(1,0)node{$5$};
\draw(4,3)node{$4$}++(1,0)node{$7$}++(1,0)node{$10$}++(1,0)node{$13$}++(0,1)node{$9$}++(-1,0)node{$6$}++(-1,0)node{$3$}++(-1,0)node{$0$}++(0,1)node{$-4$}++(1,0)node{$-1$}++(1,0)node{$2$}++(1,0)node{$5$};
\draw(6.5,.5)node{\normalsize$\calu_1$};
\draw(.5,4.5)node{\normalsize$\call_1$};
\end{tikzpicture}
\]
}

Now let's say that $x$ is a \emph{\pp} for a partition $\la$ if
\[
\call_x\subseteq\be\la\subseteq\call_x\cup\calr_x;
\]
in other words, $\be\la$ contains all elements of $\call_x$ but no elements of $\calu_x$. So a partition with \pp $x$ can be specified by saying which elements of $\calr_x$ lie in $\be\la$.

First we show that partitions having a \pp are \stm.

\begin{propn}\label{ppminl}
Suppose $\la,\mu\in\calp$ with $\cor s\la=\cor s\mu$ and $\cor t\la=\cor t\mu$, and that $x\in\bbz$ is a \pp for $\la$. Then $\card\mu\gs\card\la$, with equality \iff $x$ is also a \pp for $\mu$.
\end{propn}

\begin{pf}
For each $r\in\calr_x$, define
\begin{align*}
\dg r\la&=\card{(r+st\bbz)\cap(\calr_x\cup\calu_x)\cap\be\la}-\card{(r+st\bbz)\cap\call_x\setminus\be\la},
\\
d_r(\la)&=\sum((r+st\bbz)\cap(\calr_x\cup\calu_x)\cap\be\la)-\sum((r+st\bbz)\cap\call_x\setminus\be\la),
\end{align*}
and define $\dg r\mu$ and $d_r(\mu)$ in the same way. First we claim that
\[
\sum_{r\in\calr_x}r\dg r\la=\sum_{r\in\calr_x}r\dg r\mu.\tag*{($\dagger$)}
\]
From the definition of $\calr_x$ we get
\begin{align*}
\sum_{r\in\calr_x}r\dg r\la-\sum_{r\in\calr_x}r\dg r\mu&=\sum_{a=1}^t\sum_{b=0}^{s-1}(x+as-bt)(\dg{x+as-bt}\la-\dg{x+as-bt}\mu)\\
&=\left(\sum_{a=1}^t(x+as)\sum_{b=0}^{s-1}(\dg{x+as-bt}\la-\dg{x+as-bt}\mu)\right)\\
&\qquad-\left(\sum_{b=0}^{s-1}bt\sum_{a=1}^t(\dg{x+as-bt}\la-\dg{x+as-bt}\mu)\right).
\end{align*}
The fact that $\cor t\la=\cor t\mu$ means (using \cref{corebij}) that $\sum_{b=0}^{s-1}(\dg{x+as-bt}\la-\dg{x+as-bt}\mu)=0$ for any $a$, so that the first summand is zero. Similarly, the fact that $\cor s\la=\cor s\mu$ means that $\sum_{a=1}^t(\dg{x+as-bt}\la-\dg{x+as-bt}\mu)=0$ for any $b$, so that the second summand is also zero and ($\dagger$) is proved.

Now we consider $d_r(\la)$ and $d_r(\mu)$. It follows from \cref{basicbeta} that
\[
\card\mu-\card\la=\sum_{r\in\calr_x}(d_r(\mu)-d_r(\la))	.\tag*{($\ddagger$)}
\]
But the definitions also give
\[
d_r(\mu)\gs r\dg r\mu
\]
with equality \iff $\be\mu$ contains every element of $(r+st\bbz)\cap\call_x$ and no elements of $(r+st\bbz)\cap\calu_x$. Summing over $r$, we obtain
\[
\sum_{r\in\calr_x}d_r(\mu)\gs\sum_{r\in\calr_x}r\dg r\mu
\]
with equality \iff $\mu$ has $x$ as a \pp. Doing the same with $\la$ in place of $\mu$ and using the assumption that $x$ is a \pp for $\la$, we get
\[
\sum_{r\in\calr_x}d_r(\la)=\sum_{r\in\calr_x}r\dg r\la.
\]
So from ($\dagger$),
\[
\sum_{r\in\calr_x}d_r(\mu)\gs\sum_{r\in\calr_x}d_r(\la),
\]
and the result follows from ($\ddagger$).
\end{pf}

\begin{cory}\label{ppmins}
Suppose $\si\in\cores s$ and $\tau\in\cores t$. If there is a partition with $s$-core $\si$, $t$-core $\tau$ and a \pp $x$, then $\mns$ is precisely the set of partitions with $s$-core $\si$, $t$-core $\tau$ and $x$ as a \pp.
\end{cory}

\subsection*{The case of $(s,t)$-cores}

\cref{ppmins} is very helpful in identifying \stmps. However, there are \stmps without a \pp. For a simple example, take $(s,t)=(2,3)$, and $\la=(5,3,1^2)$. Then $\la$ is a $3$-core, so must be $(2,3)$-minimal. But it has no pinch-point, as we see from (a portion of) its $(2,3)$-diagram.
{\footnotesize
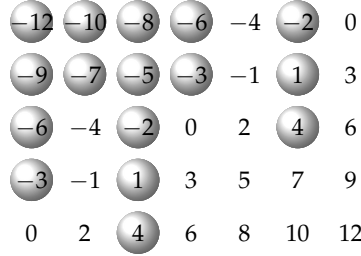
\[
\begin{tikzpicture}[scale=.7,yscale=-1]
\bkn{-1}0{-12}
\bkn00{-10}
\bkn10{-8}
\bkn20{-6}
\whn30{-4}
\bkn40{-2}
\whn500
\bkn{-1}1{-9}
\bkn01{-7}
\bkn11{-5}
\bkn21{-3}
\whn31{-1}
\bkn41{1}
\whn513
\bkn{-1}2{-6}
\whn02{-4}
\bkn12{-2}
\whn22{0}
\whn32{2}
\bkn42{4}
\whn526
\bkn{-1}3{-3}
\whn03{-1}
\bkn13{1}
\whn23{3}
\whn33{5}
\whn43{7}
\whn539
\whn{-1}4{0}
\whn04{2}
\bkn14{4}
\whn24{6}
\whn34{8}
\whn44{10}
\whn54{12}
\end{tikzpicture}
\]
}

But now we restrict to the case where $\si,\tau\in\cores{s,t}$, in which case we can show that the partitions in $\mns$ have a common \pp.

\begin{framing}{black}
For the rest of \cref{stmsec}, we assume that $\si$ and $\tau$ are $(s,t)$-cores.
\end{framing}

For any $i\in\zsz$, define
\[
\delta_i=\mfrac1s\big(\max(\be\si\cap i)-\max(\be\tau\cap i)\big).
\]
Now choose $y\in\zsz$ for which the sum
\[
\delta_0+\delta_t+\delta_{2t}+\dots+\delta_y
\]
is maximised, and let $x=\max(\be\tau\cap y)$. Call $x$ a \emph{\pvt} for $(\si,\tau)$. (Note that this depends on the order of $\si$ and $\tau$: a \pvt for $(\si,\tau)$ will not in general be a \pvt for $(\tau,\si)$.)

\begin{lemma}\label{xppst}
Suppose $\si$ and $\tau$ are $(s,t)$-cores and $x$ is a \pvt for $(\si,\tau)$. Then $x$ is a \pp for both $\si$ and $\tau$.
\end{lemma}

\begin{pf}
Let $y=x+s\bbz$. By definition $x\in\be\tau$, and the choice of $y$ means that $\delta_y\gs0$, so $x\in\be\si$ as well. Since $\si$ is an $(s,t)$-core, we get $x-as-bt\in\be\si$ for all $a,b\gs0$, so $\call_x\subseteq\be\si$. Similarly $\call_x\subseteq\be\tau$.

The definition of $x$ also means that $x+s\notin\be\tau$, and therefore $x+s+t\notin\be\tau$. The choice of $y$ implies that $\delta_{y+t}\ls0$, so $\max(\be\si\cap(y+t))\ls\max(\be\tau\cap(y+t))<x+s+t$ and therefore $x+s+t\notin\be\si$. Again using the fact that $\si$ is an $(s,t)$-core, we get $x+as+bt\notin\be\si$ for $a,b>0$, so $\be\si\cap\calu_x=\emptyset$. Similarly $\be\tau\cap\calu_x=\emptyset$, so $x$ is a \pp for $\si$ and $\tau$.
\end{pf}

In fact, $x$ is not the only common \pp of $\si$ and $\tau$: one can easily show that any two $(s,t)$-cores must have at least two \pps in common. But we will show that a \pvt for $(\si,\tau)$ is a \pp for every $\la\in\mns$. To do this, we just have to show (by \cref{ppmins}) that there is at least one partition with $s$-core $\si$, $t$-core $\tau$ and $x$ as a \pp. We do this using the theory of \zoms. First we set up a correspondence between \zoms and partitions.

\begin{lemma}\label{betapp}
Suppose $x\in\bbz$, and let $\cala_x$ denote the set of $s\times t$ \zoms with exactly $\frac12(st-s-t-1)-x$ entries equal to $1$. Then there is a bijection $\theta_x$ between $\cala_x$ and the set of partitions with \pp $x$, given by mapping a matrix $m=(m_{ij})$ to the partition with beta-set
\[
\call_x\cup\lset{x-st+it+js}{m_{ij}=1}.
\]
If $m,n\in\cala_x$, then $\theta_x(m)$ and $\theta_x(n)$ have the same $s$-core \iff $m$ and $n$ have the same row sums, while $\theta_x(m)$ and $\theta_x(n)$ have the same $t$-core \iff $m$ and $n$ have the same column sums.
\end{lemma}

\begin{pf}
If $\la\in\calp$ has $x$ as a \pp, then by \cref{basicbeta},
\[
\card{\be\la\cap\calr_x}=\card{(-\bbn)\cap(\calr_x\cup\calu_x)}-\card{\bbn_0\cap\call_x},
\]
and it is a straightforward combinatorial exercise to show that the right-hand side equals $\frac12(st-s-t-1)-x$. So $\la$ is obtained from a matrix $m\in\cala$ as described. Conversely, if $m\in\cala$, then the same calculation shows that the set
\[
B=\call_x\cup\lset{x-st+it+js}{m_{ij}=1}
\]
is a beta-set; that is, $\card{\bbn_0\cap B}=\card{(-\bbn)\setminus B}$. Clearly then the corresponding partition has $x$ as a \pp.

The statements about $s$- and $t$-cores follow from \cref{corebij}.
\end{pf}

For any partition $\la$ with $x$ as a \pp, define a composition $\rsx\la$ by setting
\[
\rsx\la_r=\card{\be\la\cap\calr_x\cap(x+rt+s\bbz)}
\]
for $r=1,\dots,s$, and $\rsx\la_r=0$ for $r>s$. Then $\rsx\la_1,\dots,\rsx\la_r$ are just the row sums of the matrix $\theta_x^{-1}(\la)$.

Now return to the situation of two $(s,t)$-cores $\si,\tau$, with $x$ being a \pvt for $(\si,\tau)$. In view of \cref{betapp}, in order to show that $x$ is a \pp of some (and hence every) partition in $\mns$, we need to show that we can find a \zom with the appropriate row and column sums.

Consider the compositions $\rsx\si$ and $\rsx\tau$. Because $\si$ is a $t$-core we have $b-t\in\be\si\cap\calr_x\cap(x+(r-1)t+s\bbz)$ for every $b\in\be\si\cap\calr_x\cap(x+rt+s\bbz)$ and every $r=2,\dots,s$, which means that $\rsx\si_1\gs\cdots\gs\rsx\si_r$, so that $\rsx\si$ is actually a partition; similarly $\rsx\tau$ is a partition, and $\card{\rsx\si}=\card{\rsx\tau}$ by \cref{betapp}. We claim that $\rsx\tau\dom\rsx\si$. If this is not true, then there is $r\in\{1,\dots,s-1\}$ for which $\si_1+\dots+\si_r>\tau_1+\dots+\tau_r$. But observe that $\rsx\si_i-\rsx\tau_i$ is the integer $\delta_{x+ti}$ defined above, so that $\delta_{x+t}+\delta_{x+2t}+\dots+\delta_{x+rt}>0$, which then contradicts the assumption that $x$ is a \pvt.

So $\rsx\si$ and $\rsx\tau$ are partitions of the same size, both with first part at most $t$ and length at most $s$, with $\rsx\tau\dom\rsx\si$. By the Gale--Ryser theorem \cite{gale,ryser}, this means that there is an $s\times t$ \zom $l$ whose row sums are $\rsx\si_1,\dots,\rsx\si_s$ and whose column sums are ${\rsx\tau}'_1,\dots,{\rsx\tau}'_t$. Now set $\la=\theta_x(l)$. Then \cref{betapp} implies that $\la$ has $x$ as a \pp, and has $s$-core $\si$ and $t$-core $\tau$.

Now \cref{ppmins} applies, and we conclude the following.

\begin{thm}\label{stcmain}
Suppose $\si,\tau\in\cores{s,t}$, and $x$ is a \pvt for $(\si,\tau)$. Then $\mns$ is the set of partitions with $s$-core $\si$, $t$-core $\tau$ and $x$ as a \pp. These partitions are in bijection with the \zoms having row sums $\rsx\si_1,\dots,\rsx\si_s$ and column sums ${\rsx\tau}'_1,\dots,{\rsx\tau}'_t$.
\end{thm}

\begin{egno}\label{maineg}
Take $s=3$, $t=4$, $\si=(1)$, $\tau=(2)$. Then
\begin{align*}
\be\si&=\{\dots,-4,-3,-2,0\},\\
\be\tau&=\{\dots,-4,-3,-2,1\}.
\end{align*}
With $\delta_i$ defined as above, we obtain
\[
\delta_0=1,\qquad\delta_1=-1,\qquad\delta_2=0,
\]
so that the unique \pvt for $(\si,\tau)$ is $-3$. The intersections of $\be\si$ and $\be\tau$ with $\calr_{-3}$ are illustrated in the following diagrams.
{\footnotesize
\[
\begin{array}{c@{\qquad\qquad\qquad}c}
\begin{tikzpicture}[scale=.6,yscale=-1]
\bkn00{-8}
\bkn10{-5}
\bkn20{-2}
\whn301
\bkn01{-4}
\whn11{-1}
\whn212
\whn315
\bkn020
\whn123
\whn226
\whn329
\end{tikzpicture}
&
\begin{tikzpicture}[scale=.6,yscale=-1]
\bkn00{-8}
\bkn10{-5}
\bkn20{-2}
\bkn301
\bkn01{-4}
\whn11{-1}
\whn212
\whn315
\whn020
\whn123
\whn226
\whn329
\end{tikzpicture}
\\
\be\si\cap\calr_{-3}
&
\be\tau\cap\calr_{-3}
\end{array}
\]
}
We see that $\rsx\si=(3,1^2)$ and $\rsx\tau=(4,1)$. So $\mns$ is in bijection with the set of $3\times4$ \zoms with row sums $3,1,1$ and column sums $2,1,1,1$. There are seven such matrices; the bijection between these matrices and the partitions $\la\in\mns$ is indicated by the following diagrams.

\[
\begin{array}{cc@{\qquad\qquad}cc}\hline
\la&\be\la\cap\calr_{-3}
&
\la&\be\la\cap\calr_{-3}
\\\hline
(10)
&
\begin{tikzpicture}[scale=.4,yscale=-1,baseline=-.5cm]
\path[use as bounding box](-.5,-.6)rectangle(3.5,2.6);
\bn00
\bn10
\bn20
\wn30
\bn01
\wn11
\wn21
\wn31
\wn02
\wn12
\wn22
\bn32
\end{tikzpicture}
&
(7,3)
&
\begin{tikzpicture}[scale=.4,yscale=-1,baseline=-.5cm]
\path[use as bounding box](-.5,-.6)rectangle(3.5,2.6);
\bn00
\bn10
\wn20
\bn30
\bn01
\wn11
\wn21
\wn31
\wn02
\wn12
\bn22
\wn32
\end{tikzpicture}
\\\hline
(6,2,1^2)
&
\begin{tikzpicture}[scale=.4,yscale=-1,baseline=-.5cm]
\path[use as bounding box](-.5,-.6)rectangle(3.5,2.6);
\bn00
\bn10
\bn20
\wn30
\wn01
\wn11
\wn21
\bn31
\bn02
\wn12
\wn22
\wn32
\end{tikzpicture}
&
(4,3,1^3)
&
\begin{tikzpicture}[scale=.4,yscale=-1,baseline=-.5cm]
\path[use as bounding box](-.5,-.6)rectangle(3.5,2.6);
\bn00
\wn10
\bn20
\bn30
\bn01
\wn11
\wn21
\wn31
\wn02
\bn12
\wn22
\wn32
\end{tikzpicture}
\\\hline
(3^3,1)
&
\begin{tikzpicture}[scale=.4,yscale=-1,baseline=-.5cm]
\path[use as bounding box](-.5,-.6)rectangle(3.5,2.6);
\bn00
\bn10
\wn20
\bn30
\wn01
\wn11
\bn21
\wn31
\bn02
\wn12
\wn22
\wn32
\end{tikzpicture}
&
(2^5)
&
\begin{tikzpicture}[scale=.4,yscale=-1,baseline=-.5cm]
\path[use as bounding box](-.5,-.6)rectangle(3.5,2.6);
\bn00
\wn10
\bn20
\bn30
\wn01
\bn11
\wn21
\wn31
\bn02
\wn12
\wn22
\wn32
\end{tikzpicture}
\\\hline
(2^2,1^6)
&
\begin{tikzpicture}[scale=.4,yscale=-1,baseline=-.5cm]
\path[use as bounding box](-.5,-.6)rectangle(3.5,2.6);
\wn00
\bn10
\bn20
\bn30
\bn01
\wn11
\wn21
\wn31
\bn02
\wn12
\wn22
\wn32
\end{tikzpicture}
\\\hline
\end{array}
\\
\]
\end{egno}

Now we determine the $t$-weight of the partitions in $\mns$, which will enable us to find their size $\nns$.

\begin{propn}\label{twtmin}
Suppose $x$ is a \pvt for $(\si,\tau)$. Then the common $t$-weight of the partitions in $\mns$ is
\[
\delta_{x+t}+2\delta_{x+2t}+\dots+(s-1)\delta_{x+(s-1)t}+s\delta_x.
\]
\end{propn}

\begin{pf}
Take $\la\in\mns$; because $\cor s\la=\si$, we get $\rsx\la=\rsx\si$. As observed above, $\delta_{x+it}=\rsx\si_i-\rsx\tau_i$ for $i=1,\dots,s$, so we need to show that the $t$-weight of $\la$ is $\sum_{i=1}^si(\rsx\la_i-\rsx\tau_i)$. In fact we show that this is true for any partition $\la$ with $t$-core $\tau$ and with $x$ as a \pp, and we do this by induction on $\wei t\la$. If $\wei t\la=0$, then $\la=\tau$ and the result is immediate. Otherwise, there is some $b\in\be\la$ with $b-t\notin\be\la$. Then both $b$ and $b-t$ lie in $\calr_x$, so $b\in x+jt+s\bbz$ for some $2\ls j\ls s$. So if we let $\mu$ denote the partition whose beta-set is obtained by replacing $b$ with $b-t$, then $\mu$ also has $t$-core $\tau$ and $x$ as a \pp. In addition $\wei t\mu=\wei t\la-1$, and
\[
\rsx\mu_i=
\begin{cases}
\rsx\la_i+1&\text{if }i=j-1\\
\rsx\la_i-1&\text{if }i=j\\
\rsx\la_i&\text{otherwise},
\end{cases}
\]
so that $\sum_{i=1}^si\rsx\mu_i=\sum_{i=1}^si\rsx\la_i-1$. So the result follows by induction.
\end{pf}

The $s$-weight of a partition $\la\in\mns$ can be determined by interchanging the roles of $s$ and $t$, or from the $t$-weight, using the fact that $\card\si+s\wei s\la=\card\la=\card\tau+t\wei t\la$.

\begin{egno}
We continue \cref{maineg}, with $-3$ being the unique \pvt for $(\si,\tau)$. According to \cref{twtmin} the $4$-weight of the partitions in $\mns$ is
\[
\delta_1+2\delta_2+3\delta_0=2.
\]
And indeed the partitions in $\mns$ have $4$-weight $2$.
\end{egno}

\subsection*{Counting partitions in $\pstn$}

Now we consider how many partitions there are in $\pstn$; in particular, we will determine exactly when $\card{\pstn}$ equals $1$ or $2$. For the case $\card{\pstn}=1$,  we have the following (which is also easily derived from the results in \cite{mfgencores}).

\begin{propn}\label{uniqueminsita}
Suppose $\si,\tau\in\cores{s,t}$, and $n\in\bbn_0$. \Tfae.
\begin{enumerate}
\item\label{uniq}
$\card{\pstn}=1$.
\item\label{sieqta}
$\si=\tau$ and $n=\card{\si}$.
\item\label{nsi}
$n=\nns=\card\si$.
\item\label{nta}
$n=\nns=\card\tau$.
\end{enumerate}
\end{propn}

\begin{pf}\indent
\begin{description}
\vspace{-\topsep}
\imps{uniq}{nsi}
If there is a unique partition of $n$ with $s$-core $\si$ and $t$-core $\tau$, then from \cref{whenpst} we certainly have $n=\nns$. Now by \cref{stcmain}, taking a \pvt $x$ for $(\si,\tau)$, there is a unique \zom with row-sum given by $\rsx\si$ and column sums given by $\rsx\tau'$. By \cite[Theorem 3.2.4]{brua}, this implies that $\rsx\si=\rsx\tau$, so that $\si=\tau$. Clearly then $\nns=\card\si$.
\imps{nsi}{sieqta}
If $\nns=\card\si$, then $\mns=\{\si\}$, because the only partition of $\card\si$ with $s$-core $\si$ is $\si$. But then $\cor t\si=\tau$, and the assumption that $\si$ is a $t$-core means that $\si=\tau$.
\imps{sieqta}{uniq}
If $\si=\tau$ and $n=\card\si$, then $\si$ is a partition of $n$ with $s$-core $\si$ and $t$-core $\tau$. Clearly it is unique, because the only partition of $\card\si$ with $s$-core $\si$ is $\si$.
\end{description}
Interchanging $s$ and $t$, we also get \ref{uniq}$\Rightarrow$\ref{nta}$\Rightarrow$\ref{sieqta}.
\end{pf}

For later use, we also want to know when there are exactly two partitions in $\pstn$. For this we need the following straightforward result on \zoms, which the author has been unable to find explicitly in the literature.

\begin{propn}\label{2zom}
Suppose $\alpha,\beta\in\calp$ with $\alpha_1,\beta_1\ls t$ and $\alpha'_1,\beta'_1\ls s$, and suppose there are exactly two $s\times t$ \zoms with rows sums $\alpha_1,\dots,\alpha_s$ and column sums $\beta'_1,\dots,\beta'_t$. Then there is some $a\in\{1,\dots,s-1\}$ for which $\alpha_a=\alpha_{a+1}$ and
\[
\beta=(\alpha_1,\dots,\alpha_{a-1},\alpha_a+1,\alpha_a-1,\alpha_{a+2},\dots).
\]
\end{propn}

(In fact the condition given on $\alpha$ and $\beta$ is sufficient as well as necessary for there to be exactly two matrices with the given row- and column-sums, but we will not need this.)

\begin{pf}
Let $m=(m_{ij})$ and $n=(n_{ij})$ be the two matrices. By Ryser's theorem \cite[Theorem 3.1]{ryser}, $m$ and $n$ differ by an \emph{interchange}; that is, there are $1\ls a<b\ls s$ and $1\ls c<d\ls t$ such that (up to switching $m$ and $n$) $m_{ac}=m_{bd}=n_{ad}=n_{bc}=1$ and $m_{ad}=m_{bc}=n_{ac}=n_{bd}=0$, while $m$ and $n$ agree in all other positions. This means that all positions except for $(a,c)$, $(a,d)$, $(b,c)$, $(b,d)$ are \emph{invariant positions}. It is easy to see that if $(i,j)$ is an invariant $1$-position (i.e.\ $m_{ij}=n_{ij}=1$) then (because $\alpha$ and $\beta$ are decreasing sequences) all positions above and/or to the left of $(i,j)$ are invariant $1$-positions, while if $(i,j)$ is an invariant $0$-position, then all positions below and/or to the right are also invariant $0$-positions. This means that $b=a+1$ and $d=c+1$.

Furthermore, $m_{ic}=m_{id}$ for all $d\neq a,b$, and $m_{aj}=m_{bj}$ for all $j\neq c,d$, since otherwise other interchanges from $m$ would be possible, so there would be more \zoms with the same row- and columns-sums as $m$. This is enough to show that $\alpha$ and $\beta$ have the desired form.
\end{pf}

Now we can characterise when there are exactly two partitions with given size, $s$-core and $t$-core.

\begin{propn}\label{extwosita}
Suppose $\si,\tau\in\cores{s,t}$, and $n\in\bbn_0$. \Tfae.
\begin{enumerate}
\item\label{ex2}
$\card{\pstn}=2$.
\item\label{addhook}
$\tau$ is obtained from $\si$ by adding or removing a rim hook, and $n=\card\si+s=\card\tau+t$.
\item\label{d11}
$n=\nns=\card\si+s=\card\tau+t$.
\end{enumerate}
\end{propn}

\needspace{3em}
\begin{pf}\indent
\begin{description}
\vspace{-\topsep}
\imps{ex2}{addhook}
If there are only two partitions of $n$ with $s$-core $\si$ and $t$-core $\tau$, then $n=\nns$ by \cref{whenpst}. Now let $x$ be a \pvt for $(\si,\tau)$, and define the partitions $\rsx\si$ and $\rsx\tau$ as above. Then by \cref{stcmain} there is a bijection from $\mns$ to the set of \zoms with row sums $\rsx\si_1,\dots,\rsx\si_s$ and column sums $\rsx\tau'_1,\dots,\rsx\tau'_t$. So if $\card\mns=2$, then by \cref{2zom} there is some $r\in\{1,\dots,s-1\}$ such that $\rsx\si_r=\rsx\si_{r+1}$ and $\rsx\tau$ is obtained from $\rsx\si$ by replacing $\rsx\si_r,\rsx\si_{r+1}$ with $\rsx\si_r+1,\rsx\si_{r+1}-1$. If we let $b=\max(\be\si\cap(x+rt+s\bbz))$, then $\be\tau$ is obtained from $\be\si$ by replacing $b+t$ with $b+s$. Hence $\tau$ is obtained from $\si$ either by removing a rim $(t-s)$-hook or by adding a rim $(s-t)$-hook. If we define a partition $\la$ by $\be\la=\be\si\setminus\{b+t\}\cup\{b+s+t\}$, then $\la\in\pst{\card\si+s}\si\tau$, so $n=\nns=\card\la=\card\si+s=\card\tau+t$.
\item[{\rm(\ref{addhook}$\Rightarrow$\ref{ex2},\ref{d11})}]
Assuming (\ref{addhook}), there is an integer $b$ such that $\be\si$ is obtained from $\be\tau$ by replacing $b+s$ with $b+t$. Since $\si$ and $\tau$ are both $(s,t)$-cores, $b$ lies in both $\be\si$ and $\be\tau$ while $b+s+t$ lies neither in $\be\si$ nor in $\be\tau$. Now define two partitions $\la$ and $\mu$ by
\[
\be\la=\be\si\setminus\{b+t\}\cup\{b+s+t\},\qquad
\be\mu=\be\si\setminus\{b\}\cup\{b+s\}.
\]
Then $\la$ and $\mu$ both have $s$-core $\si$ and $t$-core $\tau$, so are \sitam by \cref{smallwt}. Hence $\nns=\card\la=\card\si+s=\card\tau+t$, proving (\ref{d11}). Now we claim that the only partitions in $\mns$ are $\la$ and $\mu$. If $\nu\in\mns$, then $\nu$ has $s$-weight and $t$-weight both equal to $1$, so $\be\nu$ is obtained from $\be\si$ by replacing some integer $c$ with $c+s$, and is also obtained from $\be\tau$ by replacing some integer $d$ with $d+t$. The relationship between $\be\si$ and $\be\tau$ then means that $c$ can only be $b$ or $b+t$, so that $\nu=\la$ or $\mu$. So (\ref{ex2}) is proved.
\imps{d11}{addhook}
Take $\la\in\mns$. Condition (\ref{d11}) says that $\wei s\la=\wei t\la=1$, so there are $b,c\in\bbz$ such that
\[
\be\la=\be\si\setminus\{b\}\cup\{b+s\}, \qquad
\be\tau=\be\la\setminus\{c+t\}\cup\{c\}.
\]
Since $\tau$ is an $s$-core, either $b\in\be\tau$ or $b+s\notin\be\tau$. Hence either $c=b$ or $c+t=b+s$. Either way, $\be\tau$ is obtained from $\be\si$ by replacing some integer $d$ with $d+s-t$, which proves~(\ref{addhook}).
\qedhere
\end{description}
\end{pf}

\section{Conjugation}\label{conjsec}

In \cref{stmsec} we addressed the question of when a partition is determined by its $s$-core, its $t$-core and its size. Now we consider conjugation of partitions, and address the following two questions.

\begin{enumerate}
\item
When is a partition determined up to conjugation by its $s$-core, its $t$-core and its size?
\item
When is a self-conjugate partition determined by its $s$-core, its $t$-core and its size?
\end{enumerate}

In this section we answer these questions in the case of partitions whose $s$-core and $t$-core are both $(s,t)$-cores. Then in \cref{group1} we will use group actions to extend these results to all partitions.

\begin{framing}{black}
For the rest of \cref{conjsec}, we assume that $\si$ and $\tau$ are $(s,t)$-cores.
\end{framing}

\subsection*{First question on conjugation}

We begin with our first question, for which we have already done most of the work. So take $\la\in\calp$, let $\si=\cor s\la$, $\tau=\cor t\la$, and assume $\si$ and $\tau$ are both $(s,t)$-cores. We want to know whether the partitions with $s$-core $\si$, $t$-core $\tau$ and size $\card\la$ are precisely $\la$ and $\la'$. If $\la'=\la$, then we already know the answer to this question from \cref{uniqueminsita}, so we assume $\la\neq\la'$. In order to have $\cor s{\la'}=\si$ and $\cor t{\la'}=\tau$, we need $\si'=\si$ and $\tau'=\tau$ from \cref{basichook}(\ref{conjcore}). Assuming this is the case, we first need to know when there are exactly two partitions with $s$-core $\si$, $t$-core $\tau$ and size $\card\la$, and this is answered in \cref{extwosita}. We then need to know whether these two partitions form a conjugate pair or are self-conjugate. We obtain the following result.

\begin{propn}\label{uptoconjst}
Suppose $\si,\tau\in\cores{s,t}$ and $n\in\bbn_0$. \Tfae.
\begin{enumerate}
\item\label{uptocpair}
There are exactly two partitions of $n$ with $s$-core $\si$ and $t$-core $\tau$, and they form a conjugate pair.
\item\label{uptocstm}
$\si'=\si$, $\tau'=\tau$, $\si$ is obtained from $\tau$ by adding or removing a rim hook, and $n=\card\si+s=\card\tau+t$.
\item\label{uptocw}
$\si'=\si$, $\tau'=\tau$, and $n=\nns=\card\sigma+s=\card\tau+t$.
\end{enumerate}
\end{propn}

\begin{pf}\indent
\begin{description}
\vst
\imps{uptocpair}{uptocstm}
If (\ref{uptocpair}) holds, then $\si'=\si$ and $\tau'=\tau$ from \cref{basichook}(\ref{conjcore}). The rest of (\ref{uptocstm}) follows from \cref{extwosita}.
\imps{uptocstm}{uptocw}
This follows from \cref{extwosita}.
\imps{uptocw}{uptocpair}
Assuming (\ref{uptocw}), \cref{extwosita} implies that $\card{\pstn}=2$. Since $\si'=\si$ and $\tau'=\tau$, the set $\pstn$ is closed under conjugation, so the two partitions in $\pstn$ either form a conjugate pair or are both self-conjugate. But (\ref{uptocw}) says that these partitions have $s$-weight $1$, and two partitions with $s$-weight $1$ and the same $s$-core cannot both be self-conjugate by \cref{basichook}(\ref{scwt1}). So these partitions form a conjugate pair.\qedhere
\end{description}
\end{pf}

We remark that the equivalent conditions in \cref{uptoconjst} can hold only if $s-t$ is odd: if $\si$ and $\tau$ are both \sc and differ by the addition or removal of a rim hook of length $|s-t|$, then this rim hook is symmetric about the diagonal, so must contain an odd number of nodes.

\begin{egno}
Suppose $s=4$ and $t=7$. If we take $\si=(2^2)$ and $\tau=(1)$, then the conditions in \cref{uptoconjst} are satisfied when $n=8$. And indeed the smallest partitions with $4$-core $(2^2)$ and $7$-core $(1)$ are the conjugate partitions $(6,2)$ and $(2^2,1^4)$ of size $8$.
\end{egno}

\subsection*{Second question on conjugation}

Now we address our second question on conjugation: given $\si,\tau,n$, is there a unique \scp of $n$ with $s$-core $\si$ and $t$-core $\tau$? Again, this can only happen if $\si$ and $\tau$ are both \sc. But the answer to this question is considerably more complicated than for the previous question.

Recall that we define $\ol b=-1-b$ for $b\in\bbz$, and that a partition $\la$ is \sc \iff $\ol{\be\la}=\bbz\setminus\be\la$. We now introduce two more items of notation.
\begin{itemize}
\item
We write $\cst\si\tau st$ if there is an integer $b\in\be\tau\setminus\be\si$ with $\ol b-b\equiv t\ppmod{2s}$, and $\be\si=\be\tau\setminus\{b\}\cup\{\ol b\}$. This means that $\si$ is obtained from $\tau$ either by adding a rim hook whose length is congruent to $t$ modulo $2s$, or by removing a rim hook whose length is congruent to $-t$ modulo $2s$. Observe that (assuming $\si$ and $\tau$ are both \sc) this can only happen if $t$ is odd.
\item
We write $\dst\si\tau st$ if there are distinct integers $b,c\in\be\tau\setminus\be\si$ such that $\ol b-b\equiv t\ppmod{2s}$ and $\ol c-c\equiv -s\ppmod{2t}$, and $\be\si=\be\tau\setminus\{b,c\}\cup\{\ol b,\ol c\}$.  In this case we define another $(s,t)$-core $\jn\si\tau$ by $\be{\jn\si\tau}=\be\tau\setminus\{c\}\cup\{\ol c\}$.

Note that we can have $\dst\si\tau st$ only if $s$ and $t$ are both odd.
\end{itemize}

Now we can answer our question about self-conjugate partitions in the case where $\si$ and $\tau$ are both $(s,t)$-cores.

\begin{thm}\label{uniquescgeneral}
Suppose $\si,\tau\in\cores{s,t}$ and $n\in\bbn_0$. Then $\pstn$ contains a unique \scp \iff $\si'=\si$, $\tau'=\tau$ and one of the following conditions holds.
\begin{enumerate}[label=(C\arabic*)]
\item
$\si=\tau$ and $n=\card\si$.
\item
$t$ is odd, $\cst\si\tau st$ and $n=\card\tau+t$.
\item
$s$ is odd, $\cst\tau\si ts$ and $n=\card\si+s$.
\item
$s$ and $t$ are both odd, $\dst\si\tau st$ and $n=\card{\jn\si\tau}+s+t$.
\item
$s$ and $t$ are both odd, $\si=\tau=\ka_{s,t}$ and $n=\card{\ka_{s,t}}+st$.
\end{enumerate}
\end{thm}

\begin{pf}[Proof of the ``if'' part]
We show that each of the five given conditions (together with the assumption that $\si'=\si$ and $\tau'=\tau$) implies that there is a unique \scp in~$\mns$.
\begin{enumerate}[label=(C\arabic*)]
\item
In this case $\pstn=\{\si\}$, and by assumption $\si$ is \sc.
\item
If $\cst\si\tau st$, then $\be\si$ is obtained from $\be\tau$ by replacing $b$ with $\ol b=b+t-2ks$, for some integers $b,k$. So if we define a partition $\la$ by $\be\la=\be\tau\cup\{b+t-ks\}\setminus\{b-ks\}$, then $\la\in\pst{\card\tau+t}\si\tau$. In addition $\la$ is \sc, because $b+t-2ks=\ol b$ and $b+t-ks=\ol{b-ks}$. Because $\la$ has $t$-weight $1$ it is the unique \scp of $\card\tau+t$ with $t$-core $\tau$, by \cref{basichook}(\ref{scwt1}).
\item
This is the same as the previous case, with the roles of $s,t$ and of $\si,\tau$ reversed.
\item
If $\dst\si\tau st$, then $\be\si$ is obtained from $\be\tau$ by replacing two integers $b,c$ with $\ol b,\ol c$, where $\ol b-b=t-2ks$ and $\ol c-c=-s+2lt$ for some integers $k,l$. Then $c=\frac{s-1}2-lt$ and $b=\frac{-t-1}2+kt$. Let $\upsilon=\jn\si\tau$. We will show that $\nns=\card\upsilon+s+t$, and that $\mns$ contains a unique \scp.
\clam
$1\ls l\ls\frac{s-3}2$.
\prof
If $l\ls0$, then because $\tau$ is a $t$-core and $c\in\be\tau$, we get $c+2lt\in\be\tau$, and then (because $\tau$ is an $s$-core) $c+2lt-s\in\be\tau$, i.e.\ $\ol c\in\be\tau$, a contradiction. So $l\gs1$, as claimed.

To see that $l\ls\frac{s-3}2$, suppose first that $l\gs\frac{s+1}2$. Then the fact that $c\notin\be\si$ and $\si$ is a $t$-core gives $\frac{s-1}2-\frac{s+1}2t\notin\be\si$. But then the fact that $\si$ is an $s$-core yields $\frac{-t-1}2\notin\be\si$. Using the fact that $\si$ is a $t$-core again, we get $\frac{t-1}2\notin\be\si$, but this contradicts the fact that $\si'=\si$.

So we can deduce that $l\ls\frac{s-1}2$. Symmetrically, we have $k\ls\frac{t-1}2$, i.e.\ $\ol b\gs\frac{(s-1)(1-t)}2$. If in fact $l=\frac{s-1}2$, then $\be\si\not\ni c=\frac{(s-1)(1-t)}2$, and so (since $\si$ is an $s$-core)
\[
\be\si\not\ni\tfrac{(s-1)(1-t)}2+\left(\tfrac{t-1}2-k\right)=\ol b,
\]
a contradiction. So $l\ls\frac{s-3}2$.
\malc
Similarly we get $1\ls k\ls\frac{t-3}2$. Now to show that $\mns$ contains a unique \scp, we completely classify the partitions in $\mns$. Take an arbitrary partition $\mu\in\mns$, and use the set-up from \cref{stmsec} using \pvts. Recalling the integers $\delta_i$ defined there, we get
\[
\delta_b=\delta_c=-1,\qquad\delta_{b+t}=\delta_{c+2lt}=1,
\]
while all other $\delta_i$ equal $0$. So the sum
\[
\delta_0+\delta_t+\delta_{2t}\dots+\delta_y
\]
is maximised for
\[
y\ =\ b+t+s\bbz,\ b+2t+s\bbz,\ b+3t+s\bbz,\ \dots,\ c-t+s\bbz
\]
and also for
\[
y\ =\ c+2lt+s\bbz,\ c+(2l+1)t+s\bbz,\ c+(2l+2)t+s\bbz,\ \dots,\ b-t+s\bbz.
\]
This provides a large number of \pvts for $(\si,\tau)$, namely the values $\max(\be\tau\cap y)$ for the $y$ listed above. By \cref{stcmain}, each of these \pvts is a \pp for $\mu$. This means that $\be\mu\cap y=\be\si\cap y=\be\tau\cap y$ for all $y\in\zsz$ except
\[
b+s\bbz,\ b+t+s\bbz,\qquad c+s\bbz,\ c+t+s\bbz,\ \dots,\ c+2lt+s\bbz.
\]
Symmetrically, $\be\mu\cap z=\be\si\cap z=\be\tau\cap z$ for all $z\in\ztz$ except
\[
c-s+t\bbz,\ c+t\bbz,\qquad b-2ks+t\bbz,\ b-(2k-1)s+t\bbz,\ \dots,\ b+t\bbz.
\]
So $\be\mu$ agrees with $\be\si$ and $\be\tau$ except in the two disjoint sets
\begin{align*}
\cals&=\{b-2ks,b-(2k-1)s,\dots,b\}\cup\{b+t-2ks,b+t-(2k-1)s,\dots,b+t\},\\
\calt&=\{c-s,c-s+t,\dots,c-s+2lt\}\cup\{c,c+t,\dots,c+2lt\}.
\end{align*}
These sets appear in the $(s,t)$-diagram as a $2\times(2k+1)$ rectangle and a $(2l+1)\times2$ rectangle as shown in the following diagram.
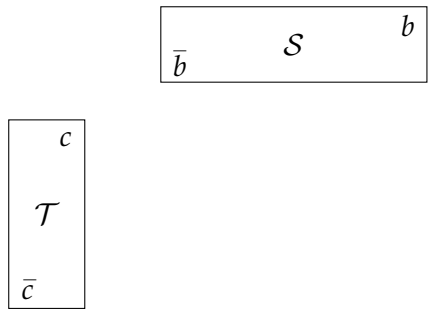
\[
\begin{tikzpicture}[scale=.5]
\draw(0,0)rectangle(2,5);
\draw(1,2.5)node{$\calt$};
\draw(0.5,0.5)node{$\ol c$};
\draw(1.5,4.5)node{$c$};
\draw(4,6)rectangle++(7,2);
\draw(7.5,7)node{$\cals$};
\draw(4.5,6.5)node{$\ol b$};
\draw(10.5,7.5)node{$b$};
\end{tikzpicture}
\]
The fact (shown above) that $1\ls l\ls\frac{s-3}2$ and $1\ls k\ls\frac{t-3}2$ means that there is no congruence class modulo $s$ or $t$ which intersects both $\cals$ and $\calt$. Given this, it is quite easy to find all $\mu\in\mns$. We consider $\be\mu\cap\cals$ first. Since $\cor s\mu=\si$, the upper row of $\cals$ contains $2k$ elements of $\be\mu$, while the lower row contains one; since $\cor t\mu=\tau$, each column of $\cals$ contains one element of $\be\mu$. So
\begin{align*}
\be\mu\cap\cals&=\{b-2kt,b+s-2kt,\dots,b\}\setminus\{b-is\}\cup\{b-is+t\}
\\
\intertext{for some $0\ls i\ls 2k$. Similarly,}
\be\mu\cap\calt&=\{c-s,c-s+t,\dots,c-s+2lt\}\setminus\{c-s+jt\}\cup\{c+jt\}
\end{align*}
for some $0\ls j\ls 2l$. So $\mns$ contains exactly $(2k+1)(2l+1)$ partitions, corresponding to the possible choices of $i$ and $j$. Only the choice $i=k$, $j=l$ gives a \scp. Call this partition $\la$; then $\be\la=\be\upsilon\setminus\{\frac{-s-1}2,\frac{-t-1}2\}\cup\{\frac{s-1}2,\frac{t-1}2\}$, so that $\nns=\card\la=\card\upsilon+s+t$.
\item
By \cref{kappaplusst}, any partition of $\card{\kappa_{s,t}}+st$ with $s$-core and $t$-core both equal to $\kappa_{s,t}$ is obtained by adding a rim $st$-hook to $\kappa_{s,t}$, i.e.\ by replacing an element $b\in\be{\kappa_{s,t}}$ with $b+st$. Since $\kappa_{s,t}$ is \sc, only the choice $b=\frac{-st-1}2$ will yield a \scp.\phantom{\qedhere}
\end{enumerate}
\end{pf}
\begin{pf}[Proof of the ``only if'' part]
Suppose $\la$ is the unique \scp of $n$ with $s$-core $\si$ and $t$-core $\tau$. Then $\si'=\cor s\la'=\cor s{\la'}=\cor s\la=\si$, and similarly $\tau'=\tau$, so we need to show that one of conditions (C1--5) holds.

Let's define a \emph{tetrad} (for $\la$) to be a quadruple $(w,x,y,z)$ of integers such that:
\begin{itemize}
\item
$w<x<z$;
\item
$w+z=x+y$;
\item
$w\equiv x\ppmod s$;
\item
$w\equiv y\ppmod t$;
\item
$\be\la\cap\{w,x,y,z\}$ equals \rlap{either $\{w,z\}$ or $\{x,y\}$.}
\end{itemize}
Let's say that this tetrad is \emph{positive} if $w,z\in\be\la$, or \emph{negative} if $x,y\in\be\la$. Observe that if $(w,x,y,z)$ is a positive tetrad, then $(\ol z,\ol y,\ol x,\ol w)$ is a negative tetrad.

\clamno\newcounter{claimcong}\setcounter{claimcong}{\value{claimcount}}
Suppose $h>0$. In any congruence class modulo $s$, there are at least as many $u$ for which $u\in\be\la\not\ni u+ht$ as there are $u$ for which $u\notin\be\la\ni u+ht$.
\prof
Since $\si$ is a $t$-core, there is no $u$ such that $u\notin\be\si\ni u+ht$, so the claim is certainly true with $\la$ replaced by $\si$. Now the claim for $\la$ follows from the fact that $\cor s\la=\si$, using \cref{corebij}.
\malc
The same statement holds with $s$ and $t$ interchanged.

\clamno\newcounter{claimtetrad}\setcounter{claimtetrad}{\value{claimcount}}
If $(w,x,y,z)$ is a positive tetrad, then $z$ equals either $\ol x$ or $\ol y$.
\prof
Suppose first that $\{w,x,y,z\}\cap\{\ol z,\ol y,\ol x,\ol w\}=\emptyset$. If $x\neq y$, then we can form a new \scp by replacing the elements $w,z,\ol x,\ol y$ in $\be\la$ with $x,y,\ol w,\ol z$. This new partition will have the same size, $s$-core and $t$-core as $\la$, contradicting our assumption that $\la$ is unique. On the other hand, if $x=y$, then $w,x,y,z$ are congruent modulo $st$ and we construct a partition $\mu$ by replacing $z,\ol x$ with $x,\ol z$ in $\be\la$. Then $\mu$ is \sc, and is obtained from $\la$ by removing an even number of $st$-hooks (so in particular has the same $s$-core and $t$-core as $\la$). Now there are several ways we can add these $st$-hooks back on to $\mu$ to create \scps, which will all have the same size, $s$-core and $t$-core as $\la$. Again, we have a contradiction.

So instead we must have $\{w,x,y,z\}\cap\{\ol z,\ol y,\ol x,\ol w\}\neq\emptyset$. Given the relationships between $w,x,y,z$ and the fact that $w,\ol x,\ol y,z\in\be\la\not\ni\ol w,x,y,\ol z$, there are four possibilities: $z=\ol x$, $z=\ol y$, $w=\ol x$ or $w=\ol y$.

Suppose $w=\ol y$. Then $\ol z,w,x$ are congruent modulo $s$, with $\ol x=\ol z+ht$, $\ol w=w+ht$ and $z=x+ht$ for some positive integer $h$. Moreover, $\ol x,w,z\in\be\la\not\ni\ol z,\ol w,x$. By Claim \arabic{claimcong} there must be some $a\neq w$ with $a\equiv w\ppmod s$ and $a\in\be\la\not\ni a+ht$. In fact (by replacing $a$ with $-a-ht$ if necessary) we can assume $a<w$. But now we have another positive tetrad $(a,x,a+ht,z)$ which is disjoint from $(\ol z,\ol{a+ht},\ol x,\ol a)$, and we get a contradiction as in the last paragraph. So $w$ cannot equal $\ol y$, and symmetrically $w$ cannot equal $\ol x$. So our claim that $z\in\{\ol x,\ol y\}$ for any positive tetrad $(w,x,y,z)$ is proved.
\malc
Correspondingly, in any negative tetrad $(w,x,y,z)$, we have $w\in\{\ol x,\ol y\}$. We split the remainder of the proof into two cases.
\begin{description}
\item[Case A: $\be\la$ does not contain $\frac{st-1}2+kst$ for any non-negative integer $k$]\indent\\
Under this assumption, we make the following claim.
\clamno
If $(w,x,y,z)$ is a positive tetrad with $z=\ol x$, then $z=\frac{t-1}2$.
\prof
To see this, first observe that because $z=\ol x$, $x\equiv z\ppmod t$ and $x<z$, we have $z=\frac{t-1}2+ht$ for some $h\gs0$. (In particular, $t$ is odd.) Assume for a contradiction that $h>0$. Because $\la$ is \sc, one of the integers $\frac{t-1}2$ and $\frac{-t-1}2$ (call it $x'$) does not lie in $\be\la$. Now we have a pair $x'<z$ with $x'\equiv z\ppmod t$ and $x'\notin\be\la\ni z$. Now Claim \arabic{claimcong} implies that $x'$ and $z$ belong to a tetrad: either a positive tetrad $(w',x',y',z)$ or a negative tetrad $(x',w',z,y')$. The assumption that $h>0$ means that $z\neq\ol{x'}$, so from Claim \arabic{claimtetrad} (and the statement immediately following it) we get either $z=\ol{y'}$ or $x'=\ol{w'}$. If $x'=\ol{w'}$, then $x'\equiv\frac{s-1}2\ppmod s$; but this is not true for either of the two possible values of $x'$. So instead $z=\ol{y'}$, which gives $z\equiv\frac{s-1}2\ppmod s$. But then $z$ has the form $\frac{st-1}2+kst$ with $k\gs0$, contrary to our current assumption.
\malc

Symmetrically, if $(w,x,y,z)$ is a positive tetrad with $z=\ol y$, then $z=\frac{s-1}2$. Combining what we have seen so far, we can deduce the following: if $u<v$ with $u\equiv v\ppmod t$ and $u\notin\be\la\ni v$, then either $u=\frac{-t-1}2$ and $v=\frac{t-1}2$, or $u=\frac{-s-1}2$, or $v=\frac{s-1}2$. The same statement holds with $s$ and $t$ interchanged throughout. So we can extract a lot of information about the sets $\be\la\cap i$, for $i\in\ztz$.
\begin{itemize}
\item
If $t$ is odd and $\frac{t-1}2\in\be\la$, then
\[
\be\la\cap(\tfrac{t-1}2+t\bbz)=\lset{\tfrac{t-1}2-ht}{h>1}\cup\{\tfrac{t-1}2\}.
\]
\item
If $s$ is odd and $\frac{s-1}2\in\be\la$, then
\begin{align*}
\be\la\cap(\tfrac{s-1}2+t\bbz)&=\lset{\tfrac{s-1}2-ht}{h>a}\cup\{\tfrac{s-1}2\},\\
\be\la\cap(\tfrac{-s-1}2+t\bbz)&=\lset{\tfrac{-s-1}2+ht}{h\ls a}\setminus\{\tfrac{-s-1}2\}
\end{align*}
for some $a>0$. (The same variable $a$ appears in both equations because $\la$ is \sc.)
\item
In all other cases, and for all other $i\in\ztz$, $\be\la\cap i$ has the form $\lset{r-ht}{h\gs0}$ for some integer $r$.
\end{itemize}
The same statement holds with $s$ and $t$ interchanged.

We can now determine which of cases (C1--5) holds, depending on whether $\frac{s-1}2$ and $\frac{t-1}2$ lie in $\be\la$.

If neither $\frac{s-1}2$ nor $\frac{t-1}2$ lies in $\be\la$, then the statements above about $\be\la$ show that $\la$ is an $(s,t)$-core, so that $\si=\tau=\la$ and (C1) holds.

If $\frac{t-1}2$ lies in $\be\la$ but $\frac{s-1}2$ does not, then certainly $t$ is odd, and the statements above show that $\wei t\la=1$, so $\la\in\mns$. $\be\tau$ is obtained from $\be\la$ by replacing $\frac{t-1}2$ with $\frac{-t-1}2$. So if we let $b$ denote the largest element of $\be\la\cap(\frac{-t-1}2+s\bbz)$, then $b=\frac{-t-1}2+as$ with $a>0$, and $\be\si$ is obtained from $\be\la$ by replacing $\frac{t-1}2$ and $b$ with $\ol b$ and $\frac{-t-1}2$. Hence $\cst\si\tau st$, so (C2) holds.

By interchanging $s$ and $t$ in the previous paragraph, we see that if $\frac{s-1}2$ lies in $\be\la$ but $\frac{t-1}2$ does not, then (C3) holds.

Finally suppose both $\frac{s-1}2$ and $\frac{t-1}2$ lie in $\be\la$. Then $s$ and $t$ are both odd. Let $b$ denote the largest element of $\be\tau\cap(\frac{-t-1}2+s\bbz)$, and let $c$ denote the largest element of $\be\tau\cap(\frac{s-1}2+t\bbz)$. Then $\be\tau$ is obtained from $\be\la$ by replacing $\frac{t-1}2$ with $\frac{-t-1}2$, $\frac{s-1}2$ with $c$, and $\ol c$ with $\frac{-s-1}2$. Similarly, $\be\si$ is obtained from $\be\la$ by replacing $\frac{s-1}2$ with $\frac{-s-1}2$, $b$ with $\frac{-t-1}2$ and $\frac{t-1}2$ with $\ol b$. Hence we have $\dst\si\tau st$, and $\card\la=\card\upsilon+s+t$. So (C4) holds.

\item[Case B: $\be\la$ contains $\frac{st-1}2+kst$ for some non-negative integer $k$]\indent\\
Take $k\gs0$ such that $\frac{st-1}2+kst$ lies in $\be\la$. Because $\la$ is \sc, the integer $\frac{-st-1}2-kst$ does not lie in $\be\la$. Let $\mu$ be the partition defined by $\be\mu=\be\la\setminus\{\frac{st-1}2+kst\}\cup\{\frac{-st-1}2-kst\}$. Then $\mu$ is \sc, and is obtained from $\la$ by removing a rim $(2k+1)st$-hook. If $k>0$, then there are several ways to add $2k+1$ rim $st$-hooks to $\mu$ to obtain a \scp of $n$ with $s$-core $\si$ and $t$-core $\tau$, contradicting the uniqueness of $\la$. So instead we must have $k=0$. We claim then $\mu$ must be the unique \scp of $n-st$ with $s$-core $\si$ and $t$-core $\tau$. If not, then let $\nu$ be another such partition. If $\frac{st-1}2\notin\be\nu$, then we can define a new \scp $\xi$ of $n$ by $\be\xi=\be\nu\cup\{\frac{st-1}2\}\setminus\{\frac{-st-1}2\}$. Then $\xi$ has $s$-core $\si$ and $t$-core $\tau$, contradicting the uniqueness of $\la$. If $\frac{st-1}2\in\be\nu$, then the partition $\pi$ of $n-2st$ defined by $\be\pi=\be\nu\setminus\{\frac{st-1}2\}\cup\{\frac{-st-1}2\}$ is \sc with $s$-core $\si$ and $t$-core $\tau$; there are several ways to add two rim $st$-hooks to $\pi$ to obtain another \scp of $n$ which contradicts the uniqueness of $\la$.

So $\mu$ is the unique \scp of $n-st$ with $s$-core $\si$ and $t$-core $\tau$, and $\be\mu$ does not contain $\frac{st-1}2+kst$ for any $k\gs0$. So from Case A applied to $\mu$, the triple $(\si,\tau,n-st)$ satisfies one of the conditions (C1--4). By the proof of the ``if'' part of the \lcnamecref{uniquescgeneral}, we know in each case exactly what $\mu$ is. We consider these cases one by one.
\begin{enumerate}[label=(C\arabic*)]
\item
In this case $\si=\tau=\mu$. If $\si=\kappa_{s,t}$, then $\si,\tau,n$ satisfy condition (C5). If not, then we claim that $\frac{-(s-1)(t-1)}2\in\be\si$. To see this, recall the set $\calx$ from the proof of \cref{kappaplusst}. This set defines $\kappa_{s,t}$, in the sense that
\[
\be{\kappa_{s,t}}=\bigcup_{x\in\calx}\{x,x-t,x-2t,\dots\}.
\]
Since $\si\neq\kappa_{s,t}$, $\be\si$ contains an integer not in $\be{\kappa_{s,t}}$, so contains $x+at$ for some $x\in\calx$ and $a>0$. Since $\si$ is a $t$-core, $\be\si$ then contains $x+t$. Since the elements of $\calx$ are congruent modulo $s$ and $\si$ is an $s$-core, $\be\si$ then contains $\min(\calx)+t=\frac{-(s-1)(t-1)}2$, as claimed. But now there is a tetrad $(\frac{-(s-1)(t-1)}2,\frac{t-1}2,\frac{s-1}2,\frac{st-1}2)$ for $\la$ which contradicts Claim \arabic{claimtetrad} above.
\item
In this case $\be\la$ is obtained from $\be\tau$ by replacing $\frac{-t-1}2$ and $\frac{-st-1}2$ with $\frac{t-1}2$ and $\frac{st-1}2$. We certainly have $\frac{-st+s-t-1}2\in\be\tau$: if $\tau=\kappa_{s,t}$, then this comes from the proof of \cref{kappaplusst}, and otherwise it comes from the fact that $\frac{-(s-1)(t-1)}2\in\be\tau$ (as shown in the case just above) and $\tau$ is a $t$-core. But now we have a tetrad $(\frac{-st+s-t-1}2,\frac{-t-1}2,\frac{s-1}2,\frac{st-1}2)$ for $\la$ which contradicts Claim \arabic{claimtetrad}.
\item
This is the same as the preceding case, with $s$ and $t$ interchanged.
\item
This is similar to (C2), but now using the tetrad $(\frac{-st-s-t-1}2,\frac{t-1}2,\frac{s-1}2,\frac{st-1}2)$.
\qedhere
\end{enumerate}
\end{description}
\end{pf}

\begin{egno}
Take $s=5$, $t=9$, $\si=(7,3^2,1^4)$, $\tau=(6,3^2,1^3)$. Then $\dst\si\tau st$, with $b=5$ and $c=-7$, giving $\upsilon=\jn\si\tau=(7^2,4^2,2^3)$. Then $\nns=\card\upsilon+14=42$, and $\la=(7^3,6^3,3)$ is the unique \scp in $\mns$. The integer $x=1$ is a \pvt for $(\si,\tau)$, and we illustrate the intersection of $\calr_1$ with the beta-set of each of the partitions $\si$, $\tau$, $\jn\si\tau$ and $\la$, with the sets $\cals$ and $\calt$ from the proof outlined.
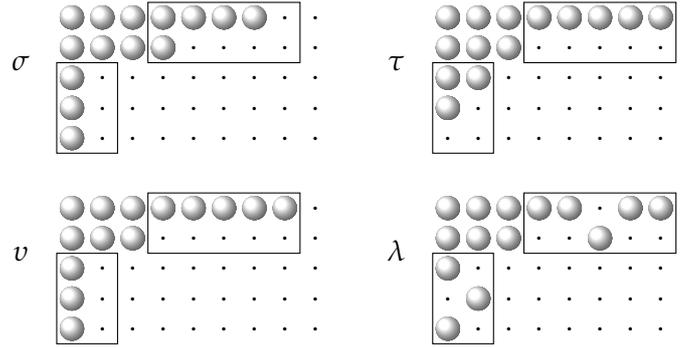
\[
\begin{array}{cc@{\qquad}cc}
\si&
\begin{tikzpicture}[scale=.4,yscale=-1,baseline=-1.1cm]
\bn11\bn12\bn13\bn14\bn15
\bn21\bn22\wn23\wn24\wn25
\bn31\bn32\wn33\wn34\wn35
\bn41\bn42\wn43\wn44\wn45
\bn51\wn52\wn53\wn54\wn55
\bn61\wn62\wn63\wn64\wn65
\bn71\wn72\wn73\wn74\wn75
\wn81\wn82\wn83\wn84\wn85
\wn91\wn92\wn93\wn94\wn95
\draw(.5,2.5)rectangle++(2,3);
\draw(3.5,.5)rectangle++(5,2);
\end{tikzpicture}
&
\tau&
\begin{tikzpicture}[scale=.4,yscale=-1,baseline=-1.1cm]
\bn11\bn12\bn13\bn14\wn15
\bn21\bn22\bn23\wn24\wn25
\bn31\bn32\wn33\wn34\wn35
\bn41\wn42\wn43\wn44\wn45
\bn51\wn52\wn53\wn54\wn55
\bn61\wn62\wn63\wn64\wn65
\bn71\wn72\wn73\wn74\wn75
\bn81\wn82\wn83\wn84\wn85
\wn91\wn92\wn93\wn94\wn95
\draw(.5,2.5)rectangle++(2,3);
\draw(3.5,.5)rectangle++(5,2);
\end{tikzpicture}
\\[42pt]
\upsilon&
\begin{tikzpicture}[scale=.4,yscale=-1,baseline=-1.1cm]
\bn11\bn12\bn13\bn14\bn15
\bn21\bn22\wn23\wn24\wn25
\bn31\bn32\wn33\wn34\wn35
\bn41\wn42\wn43\wn44\wn45
\bn51\wn52\wn53\wn54\wn55
\bn61\wn62\wn63\wn64\wn65
\bn71\wn72\wn73\wn74\wn75
\bn81\wn82\wn83\wn84\wn85
\wn91\wn92\wn93\wn94\wn95
\draw(.5,2.5)rectangle++(2,3);
\draw(3.5,.5)rectangle++(5,2);
\end{tikzpicture}
&
\la&
\begin{tikzpicture}[scale=.4,yscale=-1,baseline=-1.1cm]
\bn11\bn12\bn13\wn14\bn15
\bn21\bn22\wn23\bn24\wn25
\bn31\bn32\wn33\wn34\wn35
\bn41\wn42\wn43\wn44\wn45
\bn51\wn52\wn53\wn54\wn55
\wn61\bn62\wn63\wn64\wn65
\bn71\wn72\wn73\wn74\wn75
\bn81\wn82\wn83\wn84\wn85
\wn91\wn92\wn93\wn94\wn95
\draw(.5,2.5)rectangle++(2,3);
\draw(3.5,.5)rectangle++(5,2);
\end{tikzpicture}
\end{array}
\]
\end{egno}

\section{Group actions}\label{group1}

In this section we extend the results of \cref{stmsec,conjsec} to the case where $\si$ and $\tau$ are not necessarily $(s,t)$-cores, using actions of affine symmetric groups. These actions were introduced in \cite{mfcores}, and studied further in \cite{mfgencores,mfwtarmstrong}, and we recall the essential details here, taking the exposition from \cite{mfgencores}. Throughout, we write $\ci st$ for the integer $\frac{(s-1)(t-1)}2$.

\subsection*{The affine symmetric group}

Let $\weyl s$ denote the Coxeter group of type $\tilde A_{s-1}$. This has generators $w_i$ for $i\in\zsz$, and relations
\begin{alignat*}2
w_i^2&=1\qquad&&\text{for each }i,\\
w_iw_j&=w_jw_i&&\text{if }j\neq i\pm1,\\
w_iw_jw_i&=w_jw_iw_j&\qquad&\text{if }j=i+1\neq i-1.
\end{alignat*}

We define the \emph{level $t$ action} of $\weyl s$ on $\bbz$ by
\[
w_in=\begin{cases}
n+t&\text{if }n\in(i-1)t-\ci st\\
n-t&\text{if }n\in it-\ci st\\
n&\text{otherwise}
\end{cases}\tag*{for each $i\in\zsz$.}
\]
This naturally yields an action of $\weyl s$ on the set of all subsets of $\bbz$. This action sends beta-sets to beta-sets, so we get a level $t$ action of $\weyl s$ on the set of all partitions, defined by
\[
\be{w\la}=w\be\la
\]
for $w\in\weyl s$ and $\la\in\calp$. From now on we will refer to this simply as ``the action'' of $\weyl s$ on $\calp$, and whenever we write $w\la$ for $w\in\weyl s$ and $\la\in\calp$, we will be referring to this action. The action has several nice properties, summarised in the following \lcnamecref{actionbasic}.

\begin{lemma}[\xcite{mfgencores}{Lemma 3.1}]\label{actionbasic}
Suppose $\la\in\calp$ and $w\in\weyl s$.
\begin{enumerate}
\item
$\cor t{w\la}=\cor t\la$.
\item\label{wspres}
$\wei s{w\la}=\wei s\la$.
\item
$\cor s{w\la}=w(\cor s\la)$.
\end{enumerate}
\end{lemma}

A consequence of part \ref{wspres} of the \lcnamecref{actionbasic} is that the action of $\weyl s$ restricts to an action on $\cores s$. By interchanging $s$ and $t$ in the definitions, we get a level $s$ action of $\weyl t$ on $\calp$. This commutes with the level $t$ action of $\weyl s$, so we have an action of $\wst$.

\subsection*{\stmps}

We want to show that the set of \stmps is preserved under the action of $\weyl s$. First we observe a useful result about the sizes of partitions. For $\la\in\calp$, define $\xx\la=\card\la-\card{\cor s\la}-\card{\cor t\la}$.

\begin{propn}\label{xlapres}
Suppose $\la\in\calp$ and $w\in\wst$. Then $\xx{w\la}=\xx\la$.
\end{propn}

\begin{pf}
Assume first that $w\in\weyl s$. Then
\begin{align*}
\xx{w\la}&=\card{w\la}-\card{\cor s{w\la}}-\card{\cor t{w\la}}\\
&=s\wei s{w\la}-\card{\cor t{w\la}}\\
&=s\wei s\la-\card{\cor t\la}\tag*{by \cref{actionbasic}}\\
&=\card{\la}-\card{\cor s{\la}}-\card{\cor t{\la}}\\
&=\xx\la,
\end{align*}
as required. The case where $w\in\weyl t$ is proved in the same way, with $s$ and $t$ interchanged, and the general case follows by combining these two cases.
\end{pf}

\begin{cory}\label{stxorbit}
The set of \stmps is preserved under the action of $\wst$ on~$\calp$.
\end{cory}

\begin{pf}
To see this, take $w\in\wst$ and partitions $\la$ and $\mu$ with $w\la=\mu$. If $\la$ is not an \stmp, let $\la^-$ be a smaller partition with the same $s$-core and $t$-core as $\la$. Let $\mu^-=w\la^-$. Then we claim that $\mu^-$ has the same $s$-core and $t$-core as $\mu$, but is smaller, so that $\mu$ is not \stm. We assume first that $w\in\weyl s$. Then by \cref{actionbasic}
\[
\cor s{\mu^-}=w\cor s{\la^-}=w\cor s\la=\cor s\mu,
\]
\[
\cor t{\mu^-}=\cor t{\la^-}=\cor t\la=\cor t\mu,
\]
and so by \cref{xlapres}
\[
\card{\mu^-}=\card\mu-\card\la+\card{\la^-}<\card\mu,
\]
as claimed. The case where $w\in\weyl t$ is the same, but with $s$ and $t$ interchanged, and the general case follows by combining these two cases.
\end{pf}

We obtain the following corollary on the sizes of \stmps.

\begin{cory}\label{wnns}
Suppose $\si\in\cores s$, $\tau\in\cores t$, $u\in\weyl s$ and $v\in\weyl t$. Then
\[
\nn{u\sigma}{v\tau}=\nns-\card\si-\card\tau+\card{u\si}+\card{v\tau}.
\]
\end{cory}

\begin{pf}
Take $\la\in\mns$, and let $\mu=uv\la$. Then $\cor s\mu=u\si$ and $\cor t\mu=v\tau$ by \cref{actionbasic}, so that $\mu\in\mn{u\si}{v\tau}$, by \cref{stxorbit}. Hence
\begin{align*}
\nn{u\si}{v\tau}-\nns&=\card\mu-\card\la\\
&=\left(\xx\mu+\card{u\si}+\card{v\tau}\right)-\left(\xx\la+\card\si+\card\tau\right)\\
&=\card{u\si}+\card{v\tau}-\card\si-\card\tau\tag*{by \cref{xlapres}.\qedhere}
\end{align*}
\end{pf}

In order to use these results, we need the following \lcnamecref{coresorbit}.

\begin{lemma}\label{coresorbit}
Suppose $\calo$ is an orbit for the action of $\weyl s$ on $\cores s$. Then $\calo$ contains exactly one $(s,t)$-core.
\end{lemma}

\begin{pf}
Take $\la\in\calo$. Then $\calo$ contains $\cor t\la$ by \cite[Proposition 4.6]{mfgencores}, and by \cref{olssonthm} this is an $(s,t)$-core. By \cite[Corollary 4.7]{mfgencores} the orbit for $\wst$ containing $\calo$ contains only one $(s,t)$-core, so $\calo$ certainly contains only one.
\end{pf}

This provides a method for working out $\nn{\si}{\tau}$ for any $\si\in\cores s$ and $\tau\in\cores t$: find $u\in\weyl s$ and $v\in\weyl t$ such that $u\si$ and $v\tau$ are both $(s,t)$-cores, calculate $\nn{u\si}{v\tau}$ using \cref{twtmin}, and apply \cref{wnns}. We can find $\card{\mns}$ in the same way: by \cref{actionbasic,stxorbit}, the action of the element $uv\in\wst$ restricts to a bijection from $\mns$ to $\mn{u\si}{v\tau}$. So $\card\mns=\card{\mn{u\si}{v\tau}}$, which (modulo the difficult theory of \zoms!) can be worked out using \cref{stcmain}.

\subsection*{Orbits and stabilisers}

\cref{stxorbit} implies that $\wst$ acts on the set of \stmps, and it is interesting to ask about orbits and stabilisers for this action. First we consider the stabiliser $\stab_{\wst}(\mns)$, for given $\si,\tau$.

\begin{propn}\label{stabi}
Suppose $\si\in\cores s$ and $\tau\in\cores t$. Then $\stab_{\wst}(\mns)=\stab_{\weyl s}(\si)\times\stab_{\weyl t}(\tau)$.
\end{propn}

\begin{pf}
Take $\la\in\mns$, $u\in\weyl s$ and $v\in\weyl t$. By \cref{stxorbit}, the partition $uv\la$ lies in $\mns$ \iff it has $s$-core $\si$ and $t$-core $\tau$. By \cref{actionbasic}
\[
\cor s{uv\la}=u\si,\qquad\cor t{uv\la}=v\tau,
\]
so $uv\la\in\mns$ \iff $u\in\stab_{\weyl s}(\si)$ and $v\in\stab_{\weyl t}(\tau)$.
\end{pf}

When $\si\in\cores s$, the stabiliser $\stab_{\weyl s}(\si)$ is found in \cite[Proposition 3.7]{mfwtarmstrong}. To describe this, it is helpful to recall the \emph{$s$-set} of $\si$ introduced in \cite{mfcores}: for each $i\in\zsz$ we define $a_i$ to be the smallest element of $i$ not contained in $\be\si$, and then set $\sset\si=\lset{a_i}{i\in\zsz}$. Then \cite[Proposition 3.7]{mfwtarmstrong} says that $\card{\stab_{\weyl s}(\si)}=\prod_{i\in\ztz}\card{i\cap\sset\si}!$. If we restrict attention to the case where $\si$ is an $(s,t)$-core (which is not unduly restrictive, in view of \cref{coresorbit}), then for each $i\in\ztz$ the set $i\cap\sset\si$ is an arithmetic progression with common difference $t$, so that $\stab_{\weyl s}(\si)$ is actually a parabolic subgroup of $\weyl s$: it is generated by the elements $w_i$ for those $i\in\zsz$ satisfying $a_{it-s\circ t}=a_{(i-1)t-s\circ t}+t$. This can be seen in terms of the $(s,t)$-diagram: the generators of $\stab_{\weyl s}(\si)$ correspond to pairs of consecutive equal rows in the $(s,t)$-diagram of~$\si$. A corresponding statement holds for $\stab_{\weyl t}(\tau)$ and the columns of the $(s,t)$-diagram. Now applying the correspondence between $\mns$ and the set of \zoms described in \cref{stcmain}, we find that the action of $\stab_{\wst}(\mns)$ on $\mns$ corresponds precisely to the action on \zoms by row and column permutations.

\begin{egno}
We continue \cref{maineg}. In this case the $3$-set of $\si$ is $\{-1,1,3\}$, so that $\stab_{\weyl3}(\si)$ is a copy of the symmetric group $\sss2$, generated by $w_0$. The $4$-set of $\tau$ is $\{-1,0,2,5\}$, so that $\stab_{\weyl4}(\tau)$ is a copy of $\sss3$, generated by $w_0$ and $w_3$. The stabiliser $\stab_{\weyl s}(\si)$ then acts on $\mns$ by permuting the last two rows in each of the diagrams in \cref{maineg}, while $\stab_{\weyl4}(\tau)$ acts by permuting the last three columns. In particular, we see that the action of $\stab_{\weyl3\times\weyl4}(\mns)$ on $\mns$ has two orbits.
\end{egno}

Now we can consider the orbits of $\wst$ on the set of \stmps. For this, we can reduce to the case of $(s,t)$-cores using the following \lcnamecref{orbitreduce}.

\begin{lemma}\label{orbitreduce}
Suppose $\la$ and $\mu$ are \stmps. Then $\la$ and $\mu$ lie in the same orbit under the action of $\wst$ \iff there are $\si,\tau\in\cores{s,t}$ and $w,x\in\wst$ such that $w\la,x\mu\in\mns$ and $w\la$ and $x\mu$ lie in the same orbit under the action of $\stab_{\wst}(\mns)$.
\end{lemma}

\begin{pf}
The ``if'' part is trivial, so we prove the ``only if'' part. Let $\si=\cor t{\cor s\la}$ and $\tau=\cor s{\cor t\la}$. If $\la$ and $\mu$ lie in the same orbit under $\wst$, then $\si=\cor t{\cor s\mu}$ and $\tau=\cor s{\cor t\mu}$, by \cref{actionbasic}. In addition, $\cor s\la$ and $\si$ lie in the same orbit under the action of $\weyl s$, by \cite[Proposition 4.6]{mfgencores}, so we can find $u\in\weyl s$ such that $\si=u\cor s\la$. Similarly we can find $v\in\weyl t$ such that $v\cor t\la=\tau$. Letting $w=uv$ and applying \cref{actionbasic}, we have $\cor s{w\la}=\si$ and $\cor t{w\la}=\tau$, so that $w\la\in\mns$, by \cref{stxorbit}. In the same way we can find $x\in\wst$ with $x\mu\in\mns$. If we also take $y\in\wst$ such that $y\la=\mu$, then we claim that $xyw^{-1}\in\stab_{\wst}(\mns)$, which is all we need. Writing $xyw^{-1}$ in the form $ab$, where $a\in\weyl s$ and $b\in\weyl t$, \cref{actionbasic} gives $a\si=\si$ and $b\tau=\tau$, so that $ab\in\stab_{\weyl s}(\si)\times\stab_{\weyl t}(\tau)=\stab_{\wst}(\mns)$.
\end{pf}

As a consequence, counting the orbits of $\wst$ on the set of \stmps amounts to counting the orbits of $\stab_{\wst}(\mns)$ on $\mns$, for each pair $\si,\tau$ of $(s,t)$-cores. This in turn amounts to counting \zoms up to row- and column-equivalence. However, experiments suggest that the number of orbits of $\wst$ on the set of \stmps cannot be given by a simple formula.

\subsection*{Counting partitions in $\pstn$}

Now we use our group action to generalise the results from \cref{stmsec,conjsec}. We begin with \cref{uniqueminsita}.

\begin{thm}\label{uniqueminsitag}
Suppose $\si\in\cores s$, $\tau\in\cores t$, and $n\in\bbn_0$. \Tfae.
\begin{enumerate}
\item
$\card{\pstn}=1$.
\item
$\cor t\si=\cor s\tau$ and $n=\card\si+\card\tau-\card{\cor t\si}$.
\item
$n=\nns=\card\tau+t\wei t\si$.
\item
$n=\nns=\card\si+s\wei s\tau$.
\end{enumerate}
\end{thm}

We remark that the equivalence of (1) and (2) is shown in \cite[Section 5]{mfgencores}.

\begin{pf}
Take $u\in\weyl s$ and $v\in\weyl t$, and let $n'=n-\card\si-\card\tau+\card{u\si}+\card{v\tau}$. First we show that each of conditions 1--4 holds for the triple $(\si,\tau,n)$ \iff it holds for the triple $(u\si,v\tau,n')$.
\begin{enumerate}
\item
If $\la$ is a partition of $n$ with $s$-core $\si$ and $t$-core $\tau$, then $uv\la$ is a partition of $n'$ with $s$-core $u\si$ and $t$-core $v\tau$, by \cref{actionbasic,xlapres}. The converse is also true, so the action of $uv$ gives a bijection from $\pstn$ to $\pst{n'}{u\si}{v\tau}$.
\item
\cref{actionbasic} implies that $\cor t{u\si}=\cor t\si$ and $\cor s{v\tau}=\cor s\tau$, so $\cor t{u\la}=\cor s{v\tau}$ \iff $\cor t\si=\cor s\tau$. By the definition of $n'$ we get
\[
n'-\card{u\si}-\card{v\tau}+\card{\cor t{u\si}}=n-\card\si-\card\tau+\card{\cor t\si},
\]
so that the second condition in (2) holds for $(n,\si,\tau)$ \iff it holds for $(n',u\si,v\tau)$.
\item
\cref{wnns} implies that $n'=\nn{u\si}{v\tau}$ \iff $n=\nns$. As in part (2), $n'=\card{u\si}+\card{v\tau}-\card{\cor t{u\si}}$ \iff $n=\card\si+\card\tau-\card{\cor t\si}$, which is the same as saying $n'=\card{v\tau}+t\wei s{u\si}$ \iff $n=\card\tau+t\wei t\si$.
\item
This is similar to (3).
\end{enumerate}
By \cref{coresorbit} we can find $u\in\weyl s$ such that $u\si$ is an $(s,t)$-core. Similarly we can find $v\in\weyl t$ such that $v\tau$ is an $(s,t)$-core. By \cref{uniqueminsita} the \lcnamecref{uniqueminsitag} holds for $(u\si,v\tau,n')$, and so it holds for $(\si,\tau,n)$.
\end{pf}

\begin{rmk}
It may be more helpful to have a version of \cref{uniqueminsitag} in terms of $\la$ rather than $\si$, $\tau$ and $n$: given a partition $\la$, is it the unique partition with its size, $s$-core and $t$-core? In fact, this is answered in \cite{mfgencores}: $\la$ is unique \iff $\wei s\la=\wei s{\cor t\la}$, or equivalently $\wei t\la=\wei t{\cor s\la}$. This can also be deduced from \cref{uniqueminsitag}.
\end{rmk}

\subsection*{The hyperoctahedral group}

Now we extend the results from \cref{conjsec} involving conjugation. With conjugation of partitions in mind we need to look at a subgroup of $\wst$. The following definitions are taken from \cite[Section 4.1]{mfwtarmstrong}.

Recall that $\weyl s$ has generators $w_i$ for $i\in\zsz$. We define elements $v_a$ for $a=0,1,\dots,\lfloor s/2\rfloor$ as follows:
\[
v_a=
\begin{cases}
w_0&\text{if }a=0\\
w_aw_{-a}&\text{if }1\ls a\ls(s-2)/2\\
w_aw_{-a}w_a&\text{if }a=(s-1)/2\\
w_a&\text{if }a=s/2.
\end{cases}
\]
Let $\hyp s$ be the subgroup of $\weyl s$ generated by $v_0,\dots,v_{\lfloor s/2\rfloor}$. Then $\hyp s$ is isomorphic to the \emph{affine hyperoctahedral group}, i.e.\ the Coxeter group of type $\tilde C_{\lfloor s/2\rfloor}$. By restricting the action of $\weyl s$ on $\calp$, we obtain an action (which we also refer to simply as ``the action'') of $\hyp s$ on $\calp$. Doing the same with $s$ and $t$ interchanged, we obtain an action of $\hyp t$, and hence an action of $\hyp s\times\hyp t$, on~$\calp$.

The next result follows easily from the definition of the level $t$ action of $\weyl s$ (and explains why we need the term $-s\circ t$ in the definition of this action).

\begin{lemma}[\xcite{mfwtarmstrong}{Lemma 4.3}]\label{hypconj}
If $v\in\hyp s$ and $\la\in\calp$, then $(v\la)'=v(\la')$.
\end{lemma}

We will also need the following.

\begin{propn}[\xcite{mfwtarmstrong}{Proposition 4.7}]\label{hypcore}
Suppose $\si$ is a \sc $s$-core. Then $\cor t\si$ lies in the same orbit as $\si$ under the action of $\hyp s$.
\end{propn}

Now we can generalise \cref{uptoconjst}.

\begin{thm}\label{uptoconj}
Suppose $\si\in\cores s$, $\tau\in\cores t$ and $n\in\bbn_0$. \Tfae.
\begin{enumerate}
\item
There are exactly two partitions of $n$ with $s$-core $\si$ and $t$-core $\tau$, and they form a conjugate pair.
\item
$\si'=\si$, $\tau'=\tau$, $\cor t\si$ is obtained from $\cor s\tau$ by adding or removing a rim hook, and
\[
n=\card\si+s(\wei s\tau+1)=\card\tau+t(\wei t\si+1).
\]
\item
$\si'=\si$, $\tau'=\tau$, and
\[
n=\nns=\card\si+s(\wei s\tau+1)=\card\tau+t(\wei t\si+1).
\]
\end{enumerate}
\end{thm}

\begin{pf}
The proof strategy is the same as for \cref{uniqueminsitag}. Take $u\in\hyp s$ and $v\in\hyp t$, and let $n'=n-\card\si-\card\tau+\card{u\si}+\card{v\tau}$. We show first that each of conditions 1--3 holds for $(\si,\tau,n)$ \iff it holds for $(u\si,v\tau,n')$.
\begin{enumerate}
\item
As in the proof of \cref{uniqueminsitag}, the action of $uv$ gives a bijection from $\pstn$ to $\pst{n'}{u\si}{v\tau}$. Moreover, if a conjugate pair $\la,\la'$ of partitions belong to $\pstn$, then by \cref{hypconj} the partitions $uv\la$ and $uv(\la')$ in $\pst{n'}{u\si}{v\tau}$ form a conjugate pair.
\item
\cref{hypconj} implies that $u\si$ is \sc \iff $\si$ is, and similarly for $\tau$. \cref{actionbasic} gives $\cor t{u\si}=\cor t\si$ and $\cor s{v\tau}=\cor s\tau$, so the remaining statements in (2) are true for $\si,\tau$ \iff they are true for $u\si,v\tau$.
\item
\cref{wnns} implies that $n=\nns$ \iff $n'=\nn{u\si}{v\tau}$. The rest of (3) works in the same way as (2).
\end{enumerate}
By \cref{hypcore} we can choose $u$ so that $u\si=\cor t\si$, and in particular $u\si$ is an $(s,t)$-core. Similarly we choose $v$ such that $v\tau$ is an $(s,t)$-core. Now by \cref{uptoconjst} the \lcnamecref{uptoconj} holds for $(n',u\si,v\tau)$, so it holds for $(n,\si,\tau)$.
\end{pf}

As with \cref{uniqueminsitag}, we would like a version of \cref{uptoconj} in terms of a partition $\la$: given a non-\sc partition $\la$, is it the unique partition up to conjugation with its size, $s$-core and $t$-core?

\begin{thm}\label{uptoconjla}
Suppose $\la\in\calp$. \Tfae.
\begin{enumerate}
\item\label{onlylaconj}
$\la\neq\la'$, and the only partitions with the same size, $s$-core and $t$-core as $\la$ are $\la$ and $\la'$.
\item\label{wt1wt1}
$\cor s\la$ and $\cor t\la$ are both \sc, and
\[
\wei s\la-\wei s{\cor t\la}=\wei t\la-\wei t{\cor s\la}=1.
\]
\end{enumerate}
\end{thm}

\begin{pf}\indent
Let $n=\card\la$, $\si=\cor s\la$ and $\tau=\cor t\la$. First we prove the \lcnamecref{uptoconjla} under the assumption that $\si$ and $\tau$ are both $(s,t)$-cores.
\begin{description}
\vspace{-\topsep}
\imps{onlylaconj}{wt1wt1}
If (\ref{onlylaconj}) holds in this case, then $n,\si,\tau$ satisfy the equivalent conditions in \cref{uptoconjst}. Hence $\si'=\si$ and $\tau'=\tau$, and $n=\card\si+s=\card\tau+t$. But then
\[
\wei s\la-\wei s{\cor t\la}=\mfrac{\card\la-\card\si}s-\wei s\tau=1-0=1,
\]
and similarly for $\wei t\la-\wei t{\cor s\la}$.
\imps{wt1wt1}{onlylaconj}
If (\ref{wt1wt1}) holds, then $\la$ has $s$-weight and $t$-weight both equal to $1$, so is \stm by \cref{smallwt}. So $\card\la=\nns=\card\si+s=\card\tau+t$, so again the equivalent conditions in \cref{uptoconjst} hold.
\end{description}
The general case is derived from this special case using the action of $\hyp s\times\hyp t$, as in the proof of \cref{uptoconj}.
\end{pf}

In the same way, we can generalise \cref{uniquescgeneral}. We obtain the following.

\begin{thm}\label{uniquescgeneralg}
Suppose $\si\in\cores s$, $\tau\in\cores t$ and $n\in\bbn_0$. Then $\pstn$ contains a unique \scp \iff $\si=\si'$, $\tau=\tau'$ and one of the following occurs.
\begin{enumerate}
\item
$\cor t\si=\cor s\tau$ and $n=\card\si+\card\tau-\card{\cor t\si}$.
\item
$t$ is odd, $\cst{\cor t\si}{\cor s\tau}st$ and $n=\card\tau+t(\wei t\si+1)$.
\item
$s$ is odd, $\cst{\cor s\tau}{\cor t\si}ts$ and $n=\card\si+s(\wei s\tau+1)$.
\item
$s$ and $t$ are both odd, $\dst{\cor t\si}{\cor s\tau}st$ and
\[
n=\card{\jn{\cor t\si}{\cor s\tau}}+s(\wei s\tau+1)+t(\wei t\si+1).
\]
\item
$s$ and $t$ are both odd, $\cor t\si=\cor s\tau=\ka_{s,t}$ and $n=\card\si+\card\tau-\card{\ka_{s,t}}+st$.
\end{enumerate}
\end{thm}

\begin{egno}
Take $s=3$, $t=4$ and $n=20$. There are seven \scps of $20$. These partitions have two different $3$-cores between them, namely $(4,2,1^2)$ and $(3,1^2)$, and two different $4$-cores, namely $(4^2,2^2)$ and $(2^2)$. We consider each pair $(\si,\tau)$ in turn.
\begin{itemize}
\item
If $\si=(4,2,1^2)$ and $\tau=(4^2,2^2)$, then $\cor4\si=\varnothing=\cor3\tau$, and $20=\card\si+\card\tau-\card{\cor4\si}$, so we are in case (1) of the \lcnamecref{uniquescgeneralg}, and there is a unique \scp of $20$ with $3$-core $\si$ and $4$-core $\tau$, namely $(7,5,2^3,1^2)$.
\item
If $\si=(3,1^2)$ and $\tau=(4^2,2^2)$, then $\cor4\si=\si$ while $\cor3\tau=\varnothing$ and $\wei3\tau=4$. We can check that $\cst\varnothing{(3,1^2)}43$, and $20=\card\si+3(\wei3\tau+1)$, so we are in case (3) of the \lcnamecref{uniquescgeneralg}, and there is a unique \scp of $20$ with $3$-core $\si$ and $4$-core $\tau$, namely $(8,4,2^2,1^4)$.
\item
If $\si=(4,2,1^2)$  and $\tau=(2^2)$, then $\cor4\si=\varnothing$ and $\cor3\tau=(1)$, while $\wei s\tau=1$. So $\cor t\si\neq\cor s\tau$ and $n\neq\card\si+s(\wei s\tau+1)$, so that none of conditions 1--5 holds. And there are two \scps of $20$ with $3$-core $\si$ and $4$-core $\tau$, namely $(10,2,1^8)$ and $(5^2,4^2,2)$.
\item
If $\si=(3,1^2)$  and $\tau=(2^2)$, then $\cor4\si=\si$ and $\cor3\tau=(1)$, while $\wei s\tau=1$. Again, none of conditions 1--5 holds, and there are three \scps of $20$ with $3$-core $\si$ and $4$-core $\tau$, namely $(9,3,2,1^6)$, $(6^2,,2^4)$ and $(6,4^3,1^2)$.
\end{itemize}

\end{egno}

\end{document}